\documentclass[12pt]{article}
\usepackage{amssymb}
\usepackage{amsmath}
\usepackage{myart}

\normalbaselines
\input tcilatex
\begin{document}

\author{Yuri A. Rylov}
\title{Non-Euclidean method of the generalized geometry construction and its
application to space-time geometry.}
\date{Institute for Problems in Mechanics, Russian Academy of Sciences \\
101-1 ,Vernadskii Ave., Moscow, 117526, Russia \\
email: rylov@ipmnet.ru\\
Web site: {$http://rsfq1.physics.sunysb.edu/\symbol{126}rylov/yrylov.htm$}\\
or mirror Web site: {$http://195.208.200.111/\symbol{126}rylov/yrylov.htm$}}
\maketitle

\begin{abstract}
Non-Euclidean method of the generalized geometry construction is considered.
According to this approach any generalized geometry is obtained as a result
of deformation of the proper Euclidean geometry. The method may be applied
for construction of space-time geometries. Uniform isotropic space-time
geometry other, than that of Minkowski, is considered as an example. The
problem of the geometrical objects existence and their temporal evolution
may be considered in the constructed space-time geometry. Such a statement
of the problem is impossible in the framework of the Riemannian space-time
geometry. Existence and dynamics of microparticles is considered to be
conditioned by existence of corresponding geometrical objects and their
temporal evolution in the space-time. Geometrization of the particle mass
and its momentum is produced.
\end{abstract}

\section{Introduction}

Any (generalized) geometry is a set of propositions on properties of
geometrical objects. Geometrical object is a subset of points of the point
set $\Omega $, where the geometry is given. The number of these propositions
is very large, and labelling of these propositions by real numbers is rather
difficult, because the capacity of the set of real numbers is not sufficient
for such a labelling. Labelling by means of functions appears to be more
effective, because the set of functions is "more powerful", than the set of
real numbers. For instance, let $f$ be an integer function of integer
argument $x$. Let $f\in \lbrack 0,M-1]$ and $x\in \lbrack 0,N-1]$, where $M$
and $N$ are natural numbers. Then the number $N_{f}$ of all functions $f$ is 
$N_{f}=M^{N}$. If $f$ is a function of two integer arguments $x_{1}\in
\lbrack 0,N-1]$ and $x_{2}\in \lbrack 0,N-1]$, then the number of functions $%
N_{f}=M^{\left( N^{2}\right) }$. In other words, the number of functions $%
N_{f}$ increases much faster, then the number of values $N$ of the function
arguments and the number $M$ of the function values. Thus, labelling of the
geometry propositions by means of functions seems to be more effective, than
labelling by means of numbers.

However, investigating functions of real variables, nobody tries to
calculate the number of these functions. This calculation is used only for
Boolean functions of Boolean arguments. In this case $M=2$, $N=2$ and $%
N_{f}=2^{2}=4$. In the case of Boolean functions of two arguments, we have $%
N_{f}=2^{4}=16$. The functions of real variables are described and
investigated by subsequent combinations of simple algorithms such as:
summation, multiplication, raising to a power, taking logarithm, etc.

T-geometry is a geometry, which is constructed by means of the proper
Euclidean geometry deformation. Any geometry is a construction, which
describes the mutual disposition of points and geometrical objects on the
point set $\Omega $. The mutual disposition of points $P,Q\in \Omega $ is
described by the distance $\rho \left( P,Q\right) $ between any two points $%
P,Q\in \Omega $. In the space-time geometry the distance is real and
positive for some pairs of points, and it is imaginary for some other pairs
of points. It is more convenient and useful to use the function $\sigma
\left( P,Q\right) =\frac{1}{2}\rho ^{2}\left( P,Q\right) $, which is real
for any pair of points. The function $\sigma \left( P,Q\right) $ is known as
the world function \cite{S60}. We shall use the world function for
description of the mutual disposition of points of the point set $\Omega $.

One should expect, that giving the world function 
\begin{equation}
\sigma :\qquad \Omega \times \Omega \rightarrow \mathbb{R},\qquad \sigma
\left( P,Q\right) =\sigma \left( Q,P\right) ,\qquad \sigma \left( P,P\right)
=0,\qquad \forall P,Q\in \Omega  \label{a1.1}
\end{equation}%
for all pairs $P,Q$ of points on the point set $\Omega $, we determine the
geometry completely. Let us show this for the proper Euclidean geometry $%
\mathcal{G}_{\mathrm{E}}$, described by the Euclidean world function $\sigma
_{\mathrm{E}}$. Thereafter we expand this result to any geometry $\mathcal{G}
$, described by arbitrary world function $\sigma $, satisfying the relations
(\ref{a1.1}).

Contemporary presentation of the proper Euclidean geometry is based on the
concept of the linear vector space $V_{n}$ equipped with the scalar product
of any two vectors, given on the linear vector space. Here index $n$ means
the dimension of the linear vector space, which is defined as the maximal
number of linear independent vectors. The Euclidean $n$-dimensional point
space $E_{n}$ is obtained from the $n$-dimensional vector space $V_{n}$, if
one considers all those vectors, whose origins coincide. Then the set of all
ends of all vectors forms the Euclidean point space $\Omega =\mathbb{R}^{n}$%
. Any two points $P,Q$ form the vector $\mathbf{PQ}\equiv \overrightarrow{PQ}
$, belonging to the vector space $V_{n}=\mathbb{R}^{n}$.

The vector $\mathbf{PQ}\equiv \overrightarrow{PQ}$ is the ordered set of two
points $\left\{ P,Q\right\} ,$ $P,Q\in \mathbb{R}^{n}$. The length $%
\left\vert \mathbf{PQ}\right\vert _{\mathrm{E}}$ of the vector $\mathbf{PQ}$
is defined by the relation 
\begin{equation}
\left\vert \mathbf{PQ}\right\vert _{\mathrm{E}}^{2}=2\sigma _{\mathrm{E}%
}\left( P,Q\right)  \label{a1.2}
\end{equation}%
where index "E" means that the length of the vector is taken in the proper
Euclidean space.

The scalar product $\left( \mathbf{P}_{0}\mathbf{P}_{1}.\mathbf{P}_{0}%
\mathbf{P}_{2}\right) _{\mathrm{E}}$ of two vectors $\mathbf{P}_{0}\mathbf{P}%
_{1}$ and $\mathbf{P}_{0}\mathbf{P}_{2}$ having the common origin $P_{0}$ is
defined by the relation%
\begin{equation}
\left( \mathbf{P}_{0}\mathbf{P}_{1}.\mathbf{P}_{0}\mathbf{P}_{2}\right) _{%
\mathrm{E}}=\sigma \left( P_{0},P_{1}\right) _{\mathrm{E}}+\sigma \left(
P_{0},P_{2}\right) _{\mathrm{E}}-\sigma \left( P_{1},P_{2}\right) _{\mathrm{E%
}}  \label{a1.3}
\end{equation}%
which is obtained from the Euclidean relation%
\begin{equation}
\left\vert \mathbf{P}_{1}\mathbf{P}_{2}\right\vert _{\mathrm{E}%
}^{2}=\left\vert \mathbf{P}_{0}\mathbf{P}_{2}-\mathbf{P}_{0}\mathbf{P}%
_{1}\right\vert _{\mathrm{E}}^{2}=\left\vert \mathbf{P}_{0}\mathbf{P}%
_{2}\right\vert _{\mathrm{E}}^{2}+\left\vert \mathbf{P}_{0}\mathbf{P}%
_{1}\right\vert _{\mathrm{E}}^{2}-2\left( \mathbf{P}_{0}\mathbf{P}_{1}.%
\mathbf{P}_{0}\mathbf{P}_{2}\right) _{\mathrm{E}}  \label{a1.4}
\end{equation}%
by means of the relation (\ref{a1.2}). In particular%
\begin{equation}
\left( \mathbf{P}_{0}\mathbf{P}_{1}.\mathbf{P}_{0}\mathbf{P}_{1}\right) _{%
\mathrm{E}}=2\sigma _{\mathrm{E}}\left( P_{0},P_{1}\right) =\left\vert 
\mathbf{P}_{0}\mathbf{P}_{1}\right\vert _{\mathrm{E}}^{2}  \label{a1.4a}
\end{equation}%
Note that the relation (\ref{a1.3}) is the definition of the scalar product
via the Euclidean world function $\sigma _{\mathrm{E}}$, whereas in the
conception of the linear vector space the relation (\ref{a1.4}) is so called
cosine theorem.

In the proper Euclidean space one can define the scalar product $\left( 
\mathbf{P}_{0}\mathbf{P}_{1}.\mathbf{Q}_{0}\mathbf{Q}_{1}\right) _{\mathrm{E}%
}$ of two remote vectors $\mathbf{P}_{0}\mathbf{P}_{1}$ and $\mathbf{Q}_{0}%
\mathbf{Q}_{1}$. It is defined by the relation%
\begin{equation}
\left( \mathbf{P}_{0}\mathbf{P}_{1}.\mathbf{Q}_{0}\mathbf{Q}_{1}\right) _{%
\mathrm{E}}=\sigma \left( P_{0},Q_{1}\right) _{\mathrm{E}}+\sigma \left(
P_{1},Q_{0}\right) _{\mathrm{E}}-\sigma \left( P_{0},Q_{0}\right) _{\mathrm{E%
}}-\sigma \left( P_{1},Q_{1}\right) _{\mathrm{E}}  \label{a1.5}
\end{equation}%
which follows from evident Euclidean relation%
\begin{equation}
\left( \mathbf{P}_{0}\mathbf{P}_{1}.\mathbf{Q}_{0}\mathbf{Q}_{1}\right) _{%
\mathrm{E}}=\left( \mathbf{P}_{0}\mathbf{P}_{1}.\mathbf{P}_{0}\mathbf{Q}%
_{1}\right) _{\mathrm{E}}-\left( \mathbf{P}_{0}\mathbf{P}_{1}.\mathbf{P}_{0}%
\mathbf{Q}_{0}\right) _{\mathrm{E}}  \label{b1.5a}
\end{equation}%
and relation (\ref{a1.3}), written for two terms in rhs of (\ref{b1.5a}).%
\begin{equation*}
\left( \mathbf{P}_{0}\mathbf{P}_{1}.\mathbf{P}_{0}\mathbf{Q}_{1}\right) _{%
\mathrm{E}}=\sigma \left( P_{0},P_{1}\right) _{\mathrm{E}}+\sigma \left(
P_{0},Q_{1}\right) _{\mathrm{E}}-\sigma \left( P_{1},Q_{1}\right) _{\mathrm{E%
}}
\end{equation*}%
\begin{equation*}
\left( \mathbf{P}_{0}\mathbf{P}_{1}.\mathbf{P}_{0}\mathbf{Q}_{0}\right) _{%
\mathrm{E}}=\sigma \left( P_{0},P_{1}\right) _{\mathrm{E}}+\sigma \left(
P_{0},Q_{0}\right) _{\mathrm{E}}-\sigma \left( P_{1},Q_{0}\right) _{\mathrm{E%
}}
\end{equation*}

The necessary and sufficient condition of linear dependence of $n$ vectors ,$%
\mathbf{P}_{0}\mathbf{P}_{1}$, $\mathbf{P}_{0}\mathbf{P}_{2}$,...$\mathbf{P}%
_{0}\mathbf{P}_{n}$, defined by $n+1$ points $\mathcal{P}^{n}\equiv \left\{
P_{0},P_{1},...,P_{n}\right\} $ in the proper Euclidean space, is a
vanishing of the Gram's determinant 
\begin{equation}
F_{n}\left( \mathcal{P}^{n}\right) \equiv \det \left\vert \left\vert \left( 
\mathbf{P}_{0}\mathbf{P}_{i}.\mathbf{P}_{0}\mathbf{P}_{k}\right) _{\mathrm{E}%
}\right\vert \right\vert ,\qquad i,k=1,2,...n  \label{a1.6}
\end{equation}%
Expressing the scalar products $\left( \mathbf{P}_{0}\mathbf{P}_{i}.\mathbf{P%
}_{0}\mathbf{P}_{k}\right) _{\mathrm{E}}$ in (\ref{a1.6}) via world function 
$\sigma _{\mathrm{E}}$ by means of relation (\ref{a1.3}), we obtain
definition of linear dependence of $n$ vectors ,$\mathbf{P}_{0}\mathbf{P}%
_{1} $, $\mathbf{P}_{0}\mathbf{P}_{2}$,...$\mathbf{P}_{0}\mathbf{P}_{n}$ in
the proper Euclidean space in the form%
\begin{equation}
F_{n}\left( \mathcal{P}^{n}\right) =0  \label{a1.7}
\end{equation}%
\begin{equation}
F_{n}\left( \mathcal{P}^{n}\right) \equiv \det \left\vert \left\vert \sigma
\left( P_{0},P_{i}\right) _{\mathrm{E}}+\sigma \left( P_{0},P_{k}\right) _{%
\mathrm{E}}-\sigma \left( P_{i},P_{k}\right) _{\mathrm{E}}\right\vert
\right\vert ,\qquad i,k=1,2,...n  \label{a1.8}
\end{equation}%
Relations (\ref{a1.7}), (\ref{a1.8}) form $\sigma $-immanent definition
(i.e. the definition in terms of the world function) of the linear
dependence. This definition is obtained from the theorem on the condition of
the linear dependence of $n$ vectors in the proper Euclidean space. This
definition does not contain any reference to the linear space. It looks as a
linear dependence without a linear space.

In particular, two vectors $\mathbf{P}_{0}\mathbf{P}_{1}$, $\mathbf{Q}_{0}%
\mathbf{Q}_{1}$ are collinear (linear dependent) $\mathbf{P}_{0}\mathbf{P}%
_{1}\parallel \mathbf{Q}_{0}\mathbf{Q}_{1}$, if 
\begin{equation}
\mathbf{P}_{0}\mathbf{P}_{1}\parallel \mathbf{Q}_{0}\mathbf{Q}_{1}:\qquad
\left\vert \left\vert 
\begin{array}{cc}
\left\vert \mathbf{P}_{0}\mathbf{P}_{1}\right\vert _{\mathrm{E}}^{2} & 
\left( \mathbf{P}_{0}\mathbf{P}_{1}.\mathbf{Q}_{0}\mathbf{Q}_{1}\right) _{%
\mathrm{E}} \\ 
\left( \mathbf{Q}_{0}\mathbf{Q}_{1}.\mathbf{P}_{0}\mathbf{P}_{1}\right) _{%
\mathrm{E}} & \left\vert \mathbf{Q}_{0}\mathbf{Q}_{1}\right\vert _{\mathrm{E}%
}^{2}%
\end{array}%
\right\vert \right\vert =0  \label{a1.9}
\end{equation}%
Two vectors $\mathbf{P}_{0}\mathbf{P}_{1}$, $\mathbf{Q}_{0}\mathbf{Q}_{1}$
are in parallel, if%
\begin{equation}
\mathbf{P}_{0}\mathbf{P}_{1}\upuparrows _{\mathrm{E}}\mathbf{Q}_{0}\mathbf{Q}%
_{1}:\qquad \left( \mathbf{P}_{0}\mathbf{P}_{1}.\mathbf{Q}_{0}\mathbf{Q}%
_{1}\right) _{\mathrm{E}}=\left\vert \mathbf{P}_{0}\mathbf{P}_{1}\right\vert
_{\mathrm{E}}\cdot \left\vert \mathbf{Q}_{0}\mathbf{Q}_{1}\right\vert _{%
\mathrm{E}}  \label{a1.10}
\end{equation}

Two vectors $\mathbf{P}_{0}\mathbf{P}_{1}$, $\mathbf{Q}_{0}\mathbf{Q}_{1}$
are antiparallel, if%
\begin{equation}
\mathbf{P}_{0}\mathbf{P}_{1}\uparrow \downarrow _{\mathrm{E}}\mathbf{Q}_{0}%
\mathbf{Q}_{1}:\qquad \left( \mathbf{P}_{0}\mathbf{P}_{1}.\mathbf{Q}_{0}%
\mathbf{Q}_{1}\right) _{\mathrm{E}}=-\left\vert \mathbf{P}_{0}\mathbf{P}%
_{1}\right\vert _{\mathrm{E}}\cdot \left\vert \mathbf{Q}_{0}\mathbf{Q}%
_{1}\right\vert _{\mathrm{E}}  \label{a1.11}
\end{equation}%
Index "E" means that the scalar product and parallelism are considered in
the proper Euclidean space.

Now we can define equivalence of two vectors $\mathbf{P}_{0}\mathbf{P}_{1}$, 
$\mathbf{Q}_{0}\mathbf{Q}_{1}$ in the $\sigma $-immanent form. (i.e. in
terms of the world function $\sigma $). Two vectors $\mathbf{P}_{0}\mathbf{P}%
_{1}$, $\mathbf{Q}_{0}\mathbf{Q}_{1}$ are equivalent (equal), if 
\begin{equation}
\mathbf{P}_{0}\mathbf{P}_{1}\text{eqv}\mathbf{Q}_{0}\mathbf{Q}_{1}:\qquad
\left( \mathbf{P}_{0}\mathbf{P}_{1}\upuparrows \mathbf{Q}_{0}\mathbf{Q}%
_{1}\right) \wedge \left( \left\vert \mathbf{P}_{0}\mathbf{P}_{1}\right\vert
=\left\vert \mathbf{Q}_{0}\mathbf{Q}_{1}\right\vert \right)  \label{a1.12}
\end{equation}%
or 
\begin{eqnarray}
\mathbf{P}_{0}\mathbf{P}_{1}\text{eqv}\mathbf{Q}_{0}\mathbf{Q}_{1} &:&\qquad
\left( \left( \mathbf{P}_{0}\mathbf{P}_{1}.\mathbf{Q}_{0}\mathbf{Q}%
_{1}\right) =\left\vert \mathbf{P}_{0}\mathbf{P}_{1}\right\vert \cdot
\left\vert \mathbf{Q}_{0}\mathbf{Q}_{1}\right\vert \right)  \label{a1.12a} \\
&&\wedge \left( \left\vert \mathbf{P}_{0}\mathbf{P}_{1}\right\vert
=\left\vert \mathbf{Q}_{0}\mathbf{Q}_{1}\right\vert \right)  \label{a1.12b}
\end{eqnarray}

The property of the equivalence of two vectors in the proper Euclidean
geometry is reversible and transitive. It means 
\begin{equation}
\text{if\ \ \ }\mathbf{P}_{0}\mathbf{P}_{1}\text{eqv}\mathbf{Q}_{0}\mathbf{Q}%
_{1},\ \ \text{then \ }\mathbf{Q}_{0}\mathbf{Q}_{1}\text{eqv}\mathbf{P}_{0}%
\mathbf{P}_{1}  \label{a2.2}
\end{equation}%
\begin{equation}
\text{if\ \ \ }\left( \mathbf{P}_{0}\mathbf{P}_{1}\text{eqv}\mathbf{Q}_{0}%
\mathbf{Q}_{1}\right) \wedge \left( \mathbf{Q}_{0}\mathbf{Q}_{1}\text{eqv}%
\mathbf{R}_{0}\mathbf{R}_{1}\right) ,\ \ \text{then \ }\mathbf{P}_{0}\mathbf{%
P}_{1}\text{eqv}\mathbf{R}_{0}\mathbf{R}_{1}  \label{a2.3}
\end{equation}

However, the equivalence is reversible and transitive only in the proper
Euclidean geometry, where the property of parallelism of two vectors is
reversible and transitive. In the arbitrary generalized geometry the
property of parallelism as well as the equivalence are reversible and
intransitive, in general. Intransitivity of the equivalence property is
connected with its multivariance, when there are many vectors $\mathbf{Q}_{0}%
\mathbf{Q}_{1}$, $\mathbf{Q}_{0}\mathbf{Q}_{1}^{\prime }$, $\mathbf{Q}_{0}%
\mathbf{Q}_{1}^{\prime \prime }$,...which are equivalent to the vector $%
\mathbf{P}_{0}\mathbf{P}_{1}$, but not equivalent between themselves.
Multivariance of the equivalence property is conditioned by the fact, that
equations (\ref{a1.12a}), (\ref{a1.12b}), considered as a system of
equations for determination of the point $Q_{1}$ (at fixed points $%
P_{0},P_{1},Q_{0}$) has, many solutions, in general. It is possible also
such a situation, when equations (\ref{a1.12a}), (\ref{a1.12b}) have no
solution. Thus, the equivalence is multivariant and intransitive, in general.

In the proper Euclidean geometry the system of equations (\ref{a1.12a}), (%
\ref{a1.12b}) has always one and only one solution. In this case the
property of equivalence is single-variant, and transitive. However, already
in the Minkowski space-time geometry the equivalence of only timelike
vectors is single-variant and transitive. Equivalence of spacelike vectors
is multivariant and intransitive in the Minkowski space-time geometry.
However, nobody pays attention to this fact, because the spacelike vectors
are not used practically in applications to physics and to mechanics.

We shall distinguish between the equality relation (=) and the equivalence
relation (eqv), because the equality relation is always single-variant and
transitive, whereas the equivalence relation is multivariant and
intransitive, in general.

The sum $\mathbf{P}_{0}\mathbf{P}_{2}$ of two vectors $\mathbf{P}_{0}\mathbf{%
P}_{1}$, $\mathbf{P}_{1}\mathbf{P}_{2}$ 
\begin{equation*}
\mathbf{P}_{0}\mathbf{P}_{2}=\mathbf{P}_{0}\mathbf{P}_{1}+\mathbf{P}_{1}%
\mathbf{P}_{2}
\end{equation*}%
may be defined only in the case, when the end $P_{1}$ of the vector $\mathbf{%
P}_{0}\mathbf{P}_{1}$ coincide with the origin $P_{1}$ of the vector $%
\mathbf{P}_{1}\mathbf{P}_{2}$. However, using concept of equivalence, we may
define the sum of two arbitrary vectors $\mathbf{P}_{0}\mathbf{P}_{1}$ and $%
\mathbf{Q}_{0}\mathbf{Q}_{1}$ as a vector $\mathbf{R}_{0}\mathbf{R}_{2}$
with the origin at the point $R_{0}$ by means of the relation%
\begin{equation}
\mathbf{R}_{0}\mathbf{R}_{2}=\left( \mathbf{R}_{0}\mathbf{R}_{1}+\mathbf{R}%
_{1}\mathbf{R}_{2}\right) \text{eqv}\left( \mathbf{P}_{0}\mathbf{P}_{1}+%
\mathbf{Q}_{0}\mathbf{Q}_{1}\right)  \label{a1.16}
\end{equation}%
where vectors $\mathbf{R}_{0}\mathbf{R}_{1}$ and $\mathbf{R}_{1}\mathbf{R}%
_{2}$ are defined by the relations 
\begin{equation}
\mathbf{R}_{0}\mathbf{R}_{1}\text{eqv}\mathbf{P}_{0}\mathbf{P}_{1},\text{%
\qquad }\mathbf{R}_{1}\mathbf{R}_{2}\text{eqv}\mathbf{Q}_{0}\mathbf{Q}_{1}
\label{a1.16a}
\end{equation}
However, the sum $\mathbf{R}_{0}\mathbf{R}_{2}$ appears to be multivariant,
because of multivariance of relations (\ref{a1.16a}). Besides, the sum (\ref%
{a1.16}) depends, in general, on the order of terms in the sum, because
instead of (\ref{a1.16a}) one may use the relations%
\begin{equation}
\mathbf{R}_{0}\mathbf{R}_{1}\text{eqv}\mathbf{Q}_{0}\mathbf{Q}_{1},\text{%
\qquad }\mathbf{R}_{1}\mathbf{R}_{2}\text{eqv}\mathbf{P}_{0}\mathbf{P}_{1}
\label{a1.16b}
\end{equation}

Multiplication of the vector $\mathbf{Q}_{0}\mathbf{Q}_{1}$ by the real
number $\alpha $ is defined as follows. Vector $\mathbf{P}_{0}\mathbf{P}_{1}$
is the result of multiplication of the vector $\mathbf{Q}_{0}\mathbf{Q}_{1}$
by the real number $\alpha $, if 
\begin{equation}
\mathbf{P}_{0}\mathbf{P}_{1}\text{eqv}\left( \alpha \mathbf{Q}_{0}\mathbf{Q}%
_{1}\right) :\qquad \left\{ 
\begin{array}{l}
\left( \mathbf{P}_{0}\mathbf{P}_{1}\upuparrows \mathbf{Q}_{0}\mathbf{Q}%
_{1}\right) \wedge \left( \left\vert \mathbf{P}_{0}\mathbf{P}_{1}\right\vert
=\left\vert \alpha \right\vert \left\vert \mathbf{Q}_{0}\mathbf{Q}%
_{1}\right\vert \right) ,\ \ \text{if }\ \alpha \geq 0 \\ 
\left( \mathbf{P}_{0}\mathbf{P}_{1}\uparrow \downarrow \mathbf{Q}_{0}\mathbf{%
Q}_{1}\right) \wedge \left( \left\vert \mathbf{P}_{0}\mathbf{P}%
_{1}\right\vert =\left\vert \alpha \right\vert \left\vert \mathbf{Q}_{0}%
\mathbf{Q}_{1}\right\vert \right) ,\ \ \text{if }\ \alpha <0%
\end{array}%
\right.  \label{a1.14a}
\end{equation}%
It is possible another version of multiplication by the real number $\alpha $%
, which distinguishes from (\ref{a1.14a}) only for $\alpha =0$%
\begin{equation}
\mathbf{P}_{0}\mathbf{P}_{1}\text{eqv}\left( \alpha \mathbf{Q}_{0}\mathbf{Q}%
_{1}\right) :\qquad \left\{ 
\begin{array}{l}
\left( \mathbf{P}_{0}\mathbf{P}_{1}\upuparrows \mathbf{Q}_{0}\mathbf{Q}%
_{1}\right) \wedge \left( \left\vert \mathbf{P}_{0}\mathbf{P}_{1}\right\vert
=\left\vert \alpha \right\vert \left\vert \mathbf{Q}_{0}\mathbf{Q}%
_{1}\right\vert \right) ,\ \ \text{if }\ \alpha >0 \\ 
\mathbf{P}_{0}\mathbf{P}_{0},\ \ \text{if }\ \alpha =0 \\ 
\left( \mathbf{P}_{0}\mathbf{P}_{1}\uparrow \downarrow \mathbf{Q}_{0}\mathbf{%
Q}_{1}\right) \wedge \left( \left\vert \mathbf{P}_{0}\mathbf{P}%
_{1}\right\vert =\left\vert \alpha \right\vert \left\vert \mathbf{Q}_{0}%
\mathbf{Q}_{1}\right\vert \right) ,\ \ \text{if }\ \alpha <0%
\end{array}%
\right.  \label{a1.14}
\end{equation}%
For the proper Euclidean geometry both versions coincide.

To complete the $\sigma $-immanent description of the proper Euclidean
space, one needs to determine properties of the Euclidean world function $%
\sigma _{\mathrm{E}}$. They are presented in the second section, where one
can see, that the specific properties are different for Euclidean spaces of
different dimensions.

Presentation of the Euclidean geometry in the $\sigma $-immanent form (in
terms of the world function $\sigma _{\mathrm{E}}$) admits one to use
non-Euclidean method of the generalized geometry construction. The Euclidean
method of the geometry construction is based on derivation of the
geometrical propositions (theorems) from primordial propositions (axioms) of
the constructed geometry by means of logical reasonings and mathematical
calculations. This method is used always at the generalized geometry
construction. It reminds construction of functions of real variables by
means of simple procedures: summation, multiplication, etc.

The main defect of the Euclidean method is a necessity of a test of the
primordial axioms consistency. Such a test is a very complicated procedure,
which has been produced only for the proper Euclidean geometry. Besides, the
primordial axioms are comparatively simple only for uniform generalized
geometries. In this case the set of axioms is the same for all space
regions, described by the generalized geometry. In the case of non-uniform
geometry the set of axioms is different for different space regions.

\label{001}As far as there was only Euclidean method of the geometry
construction, a tendency appeared to prescribe the properties of the
Euclidean method to the geometry itself. As far as usually the geometry was
constructed, starting from a system of axioms, the tendency appeared to
consider any system of axioms, which contains concepts of point and of
straight line, as a kind of geometry (for instance, projective geometry,
affine geometry, etc.). In reality the Euclidean method of the geometry
construction and the system of axioms are something external with respect to
the Euclidean geometry in itself, as well as to other generalized geometries
(for instance, to the Riemannian geometry). Unfortunately, some
mathematicians could not separate the method of the geometry construction
from the geometry in itself, and this circumstance was a reason of rejection
of the geometry, constructed by the non-Euclidean method \cite{R2005b}.

We suggest a non-Euclidean method of the generalized geometry construction.
This method may be considered as a construction of the generalized geometry
by means of a deformation of the proper Euclidean geometry, when the
Euclidean world function $\sigma _{\mathrm{E}}$ is replaced by the world
function $\sigma $ of the geometry in question in all propositions of the
Euclidean geometry. As far as all propositions of the Euclidean geometry may
be labelled by the Euclidean world function $\sigma _{\mathrm{E}}$, which
describes each of such propositions completely, the replacement $\sigma _{%
\mathrm{E}}\rightarrow \sigma $ in all these propositions leads to a
construction of the generalized geometry, described by the world function $%
\sigma $.

The non-Euclidean method of the geometry construction can be carried out, if
we have the proper Euclidean geometry in the $\sigma $-immanent form. At
this method one does not need to separate primordial axioms. One does not
need to test their compatibility and to deduce other geometrical
propositions. At this approach all propositions of the proper Euclidean
geometry have equal rights. This approach reminds labelling of Boolean
functions, which appears to be possible, because of small number of the
Boolean functions. In the given case the labelling appears to be possible,
because the proper Euclidean geometry has been constructed and presented in
the $\sigma $-immanent form.

There is another analogy between the Boolean functions and the $\sigma $%
immanent presentation of the Euclidean geometry. The Boolean functions form
the mathematical tool of the formal logic. The formalism of the $\sigma $%
-immanent description forms a mathematical tool of the "geometric logic",
i.e. a system of rules for construction of any generalized geometry \cite%
{R2005a}. Maybe, on may speak about some metageometry which deals with all
possible geometry simultaneously.

Unexpected feature of the "geometric logic" is a multivariance of operations
in this logic. All operations of the conventional formal logic are
single-variant. It is convenient and customary, however, not all generalized
geometry can be constructed on the basis of single-variant logic rules.
Besides, the real space-time geometry is multivariant, and it is very
important to have a possibility of working with a multivariant "geometric
logic". Multivariance appears to be a very general property of the
generalized geometries and, in particular, of the space-time geometry.

The generalized geometry, constructed by means of a deformation of the
proper Euclidean geometry is called T-geometry (tubular geometry), because
in T-geometry straight lines are, in general, surfaces (tubes), but not
one-dimensional lines. This fact is conditioned by the multivariant
character of the parallelism in T-geometry. Indeed, the straight line $%
\mathcal{T}_{P_{0}P_{1}}$, passing through points $P_{0}$, $P_{1}$ is
defined by the relation 
\begin{equation}
\mathcal{T}_{P_{0}P_{1}}=\left\{ R|\mathbf{P}_{0}\mathbf{P}_{1}||\mathbf{P}%
_{0}\mathbf{R}\right\}  \label{a1.17}
\end{equation}%
and collinearity of vectors $\mathbf{P}_{0}\mathbf{P}_{1}$ and $\mathbf{P}%
_{0}\mathbf{R}$ is determined by one equation (\ref{a1.9})%
\begin{equation}
\left( \mathbf{P}_{0}\mathbf{P}_{1}.\mathbf{P}_{0}\mathbf{R}\right)
^{2}=\left\vert \mathbf{P}_{0}\mathbf{P}_{1}\right\vert ^{2}\cdot \left\vert 
\mathbf{P}_{0}\mathbf{R}\right\vert ^{2}  \label{a1.18}
\end{equation}%
In the $n$-dimensional space one equation (\ref{a1.18}) determines, in
general, $\left( n-1\right) $-dimensional surface. If the world function
deviates from the Euclidean world function slightly, this surface looks as a
tube.

T-geometry is interesting by its application to physics, in particular, to
the space-time geometry and dynamics. In T-geometry one can set the question
on existence of geometrical objects. In the Minkowski space-time geometry
the question on existence of a geometrical object is trivial in the sense,
that any geometrical object (any subset $\mathcal{O}$ of points of the point
set $\Omega $) may be considered as existing.

Let the set of points $\mathcal{O}_{P_{0}P_{1}...P_{n}}=\mathcal{O}\left( 
\mathcal{P}^{n}\right) $ be a geometrical object , where $\mathcal{P}%
^{n}\equiv \left\{ P_{0},P_{1},...P_{n}\right\} $ are $n+1$ characteristic
points, determining the geometrical object $\mathcal{O}\left( \mathcal{P}%
^{n}\right) $. The problem of existence of the geometrical object $\mathcal{O%
}\left( \mathcal{P}^{n}\right) $ is formulated as follows. The geometrical
object $\mathcal{O}\left( \mathcal{P}^{n}\right) $ exists at the point $%
P_{0}\in \Omega $, if at any point $Q_{0}\in \Omega $ one can construct such
a subset of points $\mathcal{O}\left( \mathcal{Q}^{n}\right) $, $\mathcal{Q}%
^{n}\equiv \left\{ Q_{0},Q_{1},...Q_{n}\right\} $, that $n\left( n+1\right)
/2$ relations take place 
\begin{equation}
\mathbf{P}_{i}\mathbf{P}_{k}\text{eqv}\mathbf{Q}_{i}\mathbf{Q}_{k},\qquad
i,k=0,1,...n,\qquad i<k  \label{a1.19}
\end{equation}%
According to (\ref{a1.12}) equivalence $\mathbf{P}_{i}\mathbf{P}_{k}$eqv$%
\mathbf{Q}_{i}\mathbf{Q}_{k}$ means two relations 
\begin{equation}
\left( \mathbf{P}_{i}\mathbf{P}_{k}.\mathbf{Q}_{i}\mathbf{Q}_{k}\right)
=\left\vert \mathbf{P}_{i}\mathbf{P}_{k}\right\vert \cdot \left\vert \mathbf{%
Q}_{i}\mathbf{Q}_{k}\right\vert ,\qquad \left\vert \mathbf{P}_{i}\mathbf{P}%
_{k}\right\vert =\left\vert \mathbf{Q}_{i}\mathbf{Q}_{k}\right\vert
\label{a1.20}
\end{equation}%
Thus, relations (\ref{a1.19}) form the system of $n(n+1)$ equations for
determination of $4n$ coordinates of points $Q_{1},Q_{2},...Q_{n}$ in the
4-dimensional space-time. Coordinates of the point $Q_{0}$ are given,
because they determine the displacement of the object $\mathcal{O}\left( 
\mathcal{Q}^{n}\right) $. In the Minkowski space-time it follows from $%
\left( \mathbf{P}_{0}\mathbf{P}_{k}\text{eqv}\mathbf{Q}_{0}\mathbf{Q}%
_{k}\right) \wedge \left( \mathbf{P}_{0}\mathbf{P}_{i}\text{eqv}\mathbf{Q}%
_{0}\mathbf{Q}_{i}\right) $ that $\mathbf{P}_{k}\mathbf{P}_{i}$eqv$\mathbf{Q}%
_{k}\mathbf{Q}_{i}$, provided all vectors are timelike. It means that not
all equations (\ref{a1.19}) are independent. Instead of $n(n+1)$ equations (%
\ref{a1.20}) we have $2n$ relations%
\begin{equation}
\left( \mathbf{P}_{0}\mathbf{P}_{k}.\mathbf{Q}_{0}\mathbf{Q}_{k}\right)
=\left\vert \mathbf{P}_{0}\mathbf{P}_{k}\right\vert \cdot \left\vert \mathbf{%
Q}_{0}\mathbf{Q}_{k}\right\vert ,\qquad \left\vert \mathbf{P}_{0}\mathbf{P}%
_{k}\right\vert =\left\vert \mathbf{Q}_{0}\mathbf{Q}_{k}\right\vert ,\qquad
k=1,2,...n  \label{a1.21}
\end{equation}%
for determination of $4n$ coordinates of points $Q_{0},Q_{1},...Q_{n}$.

The structure of the relations (\ref{a1.21}) in the Minkowski space-time is
such, that two relations 
\begin{equation}
\left( \mathbf{P}_{0}\mathbf{P}_{k}.\mathbf{Q}_{0}\mathbf{Q}_{k}\right)
=\left\vert \mathbf{P}_{0}\mathbf{P}_{k}\right\vert \cdot \left\vert \mathbf{%
Q}_{0}\mathbf{Q}_{k}\right\vert ,\qquad \left\vert \mathbf{P}_{0}\mathbf{P}%
_{k}\right\vert =\left\vert \mathbf{Q}_{0}\mathbf{Q}_{k}\right\vert
\label{a1.22}
\end{equation}%
determine uniquely four coordinates of the point $Q_{k}$, provided the
vector $\mathbf{P}_{0}\mathbf{P}_{k}$ is timelike, i.e. $\left\vert \mathbf{P%
}_{0}\mathbf{P}_{k}\right\vert ^{2}>0$.

In other uniform isotropic space-times the structure of relations (\ref%
{a1.22}) has another character. In this case two relations (\ref{a1.22}) do
not determine uniquely four coordinates of the point $Q_{k}$. Besides, the
relation $\mathbf{P}_{k}\mathbf{P}_{i}$eqv$\mathbf{Q}_{k}\mathbf{Q}_{i}$ is
not a corollary of $\left( \mathbf{P}_{0}\mathbf{P}_{k}\text{eqv}\mathbf{Q}%
_{0}\mathbf{Q}_{k}\right) \wedge \left( \mathbf{P}_{0}\mathbf{P}_{i}\text{eqv%
}\mathbf{Q}_{0}\mathbf{Q}_{i}\right) $ and relations (\ref{a1.19}) form $%
n\left( n+1\right) $ relations which are independent, in general. For two
characteristic points $Q_{0},Q_{1}$ we have $n=1$ and the number of
equations $n\left( n+1\right) =2$ is less, than the number of coordinates $%
4n=4$ of point $Q_{1}$.

In the case of three characteristic points $Q_{0},Q_{1},Q_{2}$ we have $n=2$%
, and the number of equations $n\left( n+1\right) =6$ is less, than the
number of coordinates $4n=8$ of points $Q_{1},Q_{2}$

In the case of four characteristic points $Q_{0},Q_{1},Q_{2},Q_{3}$ we have $%
n=3,$ and the number of equations $n\left( n+1\right) =12$ is equal to the
number of coordinates $4n=12$ of points $Q_{1},Q_{2},Q_{3}$.

Finally, in the case of five characteristic points $%
Q_{0},Q_{1},Q_{2},Q_{3},Q_{4}$ we have $n=4$, and the number of equations $%
n\left( n+1\right) =20$ is more than the number of coordinates $4n=16$ of
points $Q_{1},Q_{2},Q_{3},Q_{4}$.

It means, that geometrical objects, having more, than four characteristic
points do not exist in the multivariant space-time, in general.

In the second section one presents specific properties of the Euclidean
world function, which form the necessary and sufficient conditions of the
Euclideaness. The third and fourth sections are devoted to consideration of
the timelike vectors equivalence. In the fifth section the equivalence of
null vectors is considered. Construction of geometrical objects is
considered in the sixth section. Temporal evolution of the timelike straight
line segment is considered in the seventh section.

\section{Specific properties of the $n$-dimensional proper Euclidean space}

There are four conditions which are necessary and sufficient conditions of
the fact, that the world function $\sigma $ is the world function of $n$%
-dimensional Euclidean space \cite{R02}. They have the form:

I. Definition of the dimension: 
\begin{equation}
\exists \mathcal{P}^{n}\equiv \left\{ P_{0},P_{1},...P_{n}\right\} \subset
\Omega ,\qquad F_{n}\left( \mathcal{P}^{n}\right) \neq 0,\qquad F_{k}\left( {%
\Omega }^{k+1}\right) =0,\qquad k>n  \label{g2.5}
\end{equation}%
where $F_{n}\left( \mathcal{P}^{n}\right) $\ is the Gram's determinant (\ref%
{a1.6}). Vectors $\mathbf{P}_{0}\mathbf{P}_{i}$, $\;i=1,2,...n$\ are basic
vectors of the rectilinear coordinate system $K_{n}$\ with the origin at the
point $P_{0}$. The metric tensors $g_{ik}\left( \mathcal{P}^{n}\right) $, $%
g^{ik}\left( \mathcal{P}^{n}\right) $, \ $i,k=1,2,...n$\ in $K_{n}$\ are
defined by the relations 
\begin{equation}
\sum\limits_{k=1}^{k=n}g^{ik}\left( \mathcal{P}^{n}\right) g_{lk}\left( 
\mathcal{P}^{n}\right) =\delta _{l}^{i},\qquad g_{il}\left( \mathcal{P}%
^{n}\right) =\left( \mathbf{P}_{0}\mathbf{P}_{i}.\mathbf{P}_{0}\mathbf{P}%
_{l}\right) ,\qquad i,l=1,2,...n  \label{a1.5b}
\end{equation}%
\begin{equation}
F_{n}\left( \mathcal{P}^{n}\right) =\det \left\vert \left\vert g_{ik}\left( 
\mathcal{P}^{n}\right) \right\vert \right\vert \neq 0,\qquad i,k=1,2,...n
\label{g2.6}
\end{equation}

II. Linear structure of the Euclidean space: 
\begin{equation}
\sigma \left( P,Q\right) =\frac{1}{2}\sum\limits_{i,k=1}^{i,k=n}g^{ik}\left( 
\mathcal{P}^{n}\right) \left( x_{i}\left( P\right) -x_{i}\left( Q\right)
\right) \left( x_{k}\left( P\right) -x_{k}\left( Q\right) \right) ,\qquad
\forall P,Q\in \Omega  \label{a1.5a}
\end{equation}%
where coordinates $x_{i}\left( P\right) ,$\ $x_{i}\left( Q\right) ,$ $%
i=1,2,...n$\ of the points $P$ and $Q$\ are covariant coordinates of the
vectors $\mathbf{P}_{0}\mathbf{P}$, $\mathbf{P}_{0}\mathbf{Q}$ respectively,
defined by the relation 
\begin{equation}
x_{i}\left( P\right) =\left( \mathbf{P}_{0}\mathbf{P}_{i}.\mathbf{P}_{0}%
\mathbf{P}\right) ,\qquad i=1,2,...n  \label{b.12}
\end{equation}

III: The metric tensor matrix $g_{lk}\left( \mathcal{P}^{n}\right) $\ has
only positive eigenvalues 
\begin{equation}
g_{k}>0,\qquad k=1,2,...,n  \label{a1.5c}
\end{equation}

IV. The continuity condition: the system of equations 
\begin{equation}
\left( \mathbf{P}_{0}\mathbf{P}_{i}.\mathbf{P}_{0}\mathbf{P}\right)
=y_{i}\in \mathbb{R},\qquad i=1,2,...n  \label{b1.4}
\end{equation}%
considered to be equations for determination of the point $P$\ as a function
of coordinates $y=\left\{ y_{i}\right\} $,\ \ $i=1,2,...n$\ has always one
and only one solution. Conditions I -- IV contain a reference to the
dimension $n$\ of the Euclidean space.

\section{Equivalence of two vectors}

The property of the two vectors equality may be introduced in any T-geometry
by means of the relation (\ref{a1.12}). But in the arbitrary T-geometry the
equality of two vectors is intransitive, in general, because of the
parallelism multivariance. Intransitivity and multivariance of the two
vectors equality is very inconvenient in applications. We shall use the term
"equivalence" instead of the term "equality".

\textit{Definition. }Two vectors $\mathbf{P}_{0}\mathbf{P}_{1}$ and $\mathbf{%
Q}_{0}\mathbf{Q}_{1}$ are equivalent ($\mathbf{P}_{0}\mathbf{P}_{1}$eqv$%
\mathbf{Q}_{0}\mathbf{Q}_{1}$), if the conditions (\ref{a1.12}) take place%
\begin{equation}
\mathbf{P}_{0}\mathbf{P}_{1}\text{eqv}\mathbf{Q}_{0}\mathbf{Q}_{1}:\qquad
\left( \left( \mathbf{P}_{0}\mathbf{P}_{1}.\mathbf{Q}_{0}\mathbf{Q}%
_{1}\right) =\left\vert \mathbf{P}_{0}\mathbf{P}_{1}\right\vert \cdot
\left\vert \mathbf{Q}_{0}\mathbf{Q}_{1}\right\vert \right) \wedge \left(
\left\vert \mathbf{P}_{0}\mathbf{P}_{1}\right\vert =\left\vert \mathbf{Q}_{0}%
\mathbf{Q}_{1}\right\vert \right)  \label{a2.4}
\end{equation}%
It follows from (\ref{a2.4}), that if $\left( \mathbf{P}_{0}\mathbf{P}%
_{1}\right) $eqv$\left( \mathbf{Q}_{0}\mathbf{Q}_{1}\right) $, then $\left( 
\mathbf{Q}_{0}\mathbf{Q}_{1}\right) $eqv$\left( \mathbf{P}_{0}\mathbf{P}%
_{1}\right) $

\textit{Remark. }We distinguish between the equality (=) of vectors and
equivalence (eqv) of vectors. For instance, the equality $\mathbf{P}_{0}%
\mathbf{P}_{1}=\mathbf{P}_{0}\mathbf{Q}_{1}$ means, that the points $P_{1}$
and $Q_{1}$ coincide ($P_{1}=Q_{1}$). Equality $\mathbf{P}_{0}\mathbf{P}_{1}=%
\mathbf{Q}_{0}\mathbf{Q}_{1}$ means, that the point $P_{1}$ coincides with $%
Q_{1}$ ($P_{1}=Q_{1}$) and the point $P_{0}$ coincides with $Q_{0}$ ($%
P_{0}=Q_{0}$), whereas equivalence $\mathbf{P}_{0}\mathbf{P}_{1}$eqv$\mathbf{%
Q}_{0}\mathbf{Q}_{1}$ means the fulfilment of relations (\ref{a2.4}). The
point $P_{0}$ may not coincide with $Q_{0}$ and the point $P_{1}$ may not
coincide with $Q_{1}$, i.e. equalities $P_{0}=Q_{0}$ and $P_{1}=Q_{1}$ may
not take place.

The shift vector $\mathbf{P}_{0}\mathbf{Q}_{0}$ describes the shift of the
origin $P_{0}$ of the vector $\mathbf{P}_{0}\mathbf{P}_{1}$. The shift
vector $\mathbf{P}_{1}\mathbf{Q}_{1}$ describes the shift of the end $P_{1}$
of the vector $\mathbf{P}_{0}\mathbf{P}_{1}$. In the proper Euclidean space
equivalence of shift vectors $\mathbf{P}_{0}\mathbf{Q}_{0}$ and $\mathbf{P}%
_{1}\mathbf{Q}_{1}$ leads to equivalence of vectors $\mathbf{P}_{0}\mathbf{P}%
_{1}$ and $\mathbf{Q}_{0}\mathbf{Q}_{1}$ and vice versa equivalence of
vectors $\mathbf{P}_{0}\mathbf{P}_{1}$ and $\mathbf{Q}_{0}\mathbf{Q}_{1}$
leads to equivalence of their shift vectors $\mathbf{P}_{0}\mathbf{Q}_{0}$eqv%
$\mathbf{P}_{1}\mathbf{Q}_{1}$. In the general T-geometry the equivalence of
shift vectors $\mathbf{P}_{0}\mathbf{Q}_{0}$ and $\mathbf{P}_{1}\mathbf{Q}%
_{1}$ is not sufficient for equivalence of vectors $\mathbf{P}_{0}\mathbf{P}%
_{1}$ and $\mathbf{Q}_{0}\mathbf{Q}_{1}$. It is necessary once more
constraint $\left\vert \mathbf{P}_{0}\mathbf{P}_{1}\right\vert =\left\vert 
\mathbf{Q}_{0}\mathbf{Q}_{1}\right\vert $ or $\mathbf{P}_{0}\mathbf{P}%
_{1}\upuparrows \mathbf{Q}_{0}\mathbf{Q}_{1}$, to provide their equivalence.

\textit{Theorem}. Vectors $\mathbf{P}_{0}\mathbf{P}_{1}$ and $\mathbf{Q}_{0}%
\mathbf{Q}_{1}$ are equivalent, if shift vectors $\mathbf{P}_{0}\mathbf{Q}%
_{0}$ and $\mathbf{P}_{1}\mathbf{Q}_{1}$ are equivalent and $\left\vert 
\mathbf{P}_{0}\mathbf{P}_{1}\right\vert =\left\vert \mathbf{Q}_{0}\mathbf{Q}%
_{1}\right\vert $, or $\mathbf{P}_{0}\mathbf{P}_{1}\uparrow \uparrow \mathbf{%
Q}_{0}\mathbf{Q}_{1}$

Let $\left( \mathbf{P}_{0}\mathbf{Q}_{0}\text{eqv}\mathbf{P}_{1}\mathbf{Q}%
_{1}\right) $. Equivalence of $\mathbf{P}_{0}\mathbf{Q}_{0}$ and $\mathbf{P}%
_{1}\mathbf{Q}_{1}$ is written in the form of two relations 
\begin{equation}
\sigma \left( P_{0},Q_{1}\right) +\sigma \left( P_{1},Q_{0}\right) -\sigma
\left( P_{0},P_{1}\right) -\sigma \left( Q_{0},Q_{1}\right) =2\sqrt{\sigma
\left( P_{0},Q_{0}\right) \sigma \left( P_{1},Q_{1}\right) }  \label{b2.5}
\end{equation}%
\begin{equation}
\sigma \left( P_{0},Q_{0}\right) =\sigma \left( P_{1},Q_{1}\right)
\label{b2.6}
\end{equation}%
In force of (\ref{b2.6}) equation (\ref{b2.5}) may be written in the form%
\begin{equation}
\sigma \left( P_{0},Q_{1}\right) +\sigma \left( P_{1},Q_{0}\right) -\sigma
\left( P_{0},P_{1}\right) -\sigma \left( Q_{0},Q_{1}\right) =\sigma \left(
P_{0},Q_{0}\right) +\sigma \left( P_{1},Q_{1}\right)  \label{a2.6}
\end{equation}%
The relation $\mathbf{P}_{0}\mathbf{P}_{1}$eqv$\mathbf{Q}_{0}\mathbf{Q}_{1}$
is written in the form%
\begin{eqnarray}
\sigma \left( P_{0},Q_{1}\right) +\sigma \left( P_{1},Q_{0}\right) -\sigma
\left( P_{0},Q_{0}\right) -\sigma \left( P_{1},Q_{1}\right) &=&2\sqrt{\sigma
\left( P_{0},P_{1}\right) \sigma \left( Q_{0},Q_{1}\right) }  \label{a2.7} \\
\sigma \left( P_{0},P_{1}\right) -\sigma \left( Q_{0},Q_{1}\right) &=&0
\label{a2.8}
\end{eqnarray}

The difference of (\ref{a2.6}) and (\ref{a2.7}) has the form%
\begin{equation}
\sigma \left( P_{0},P_{1}\right) +\sigma \left( Q_{0},Q_{1}\right) =2\sqrt{%
\sigma \left( P_{0},P_{1}\right) \sigma \left( Q_{0},Q_{1}\right) }
\label{a2.9}
\end{equation}%
which can be reduced to the form%
\begin{equation}
\left( \sqrt{\sigma \left( P_{0},P_{1}\right) }-\sqrt{\sigma \left(
Q_{0},Q_{1}\right) }\right) ^{2}=0  \label{a2.10}
\end{equation}

Let now $\left\vert \mathbf{P}_{0}\mathbf{P}_{1}\right\vert =\left\vert 
\mathbf{Q}_{0}\mathbf{Q}_{1}\right\vert $, and relations (\ref{a2.8}), (\ref%
{a2.9}) take place. As far as the relations (\ref{b2.5}), (\ref{b2.6}) are
supposed to be fulfilled, the relation (\ref{a2.6}) is fulfilled also. The
relation (\ref{a2.7}) takes place also, because it is a sum of relations (%
\ref{a2.6}) and (\ref{a2.9}). Thus, equations (\ref{a2.7}) and (\ref{a2.8})
are fulfilled. It means that $\mathbf{P}_{0}\mathbf{P}_{1}$eqv$\mathbf{Q}_{0}%
\mathbf{Q}_{1}$.

Let now $\left\vert \mathbf{P}_{0}\mathbf{P}_{1}\right\vert \upuparrows
\left\vert \mathbf{Q}_{0}\mathbf{Q}_{1}\right\vert $, and relation (\ref%
{a2.7}) is fulfilled. As far as the relations (\ref{b2.5}), (\ref{b2.6}), (%
\ref{a2.6}) are supposed to be fulfilled, the relation (\ref{a2.9}) takes
place also, because the relation (\ref{a2.9}) is a difference of equations (%
\ref{a2.6}) and (\ref{a2.7}). Equation (\ref{a2.8}) is a corollary of (\ref%
{a2.9}). Thus, equations (\ref{a2.7}), (\ref{a2.8}) are fulfilled and $%
\mathbf{P}_{0}\mathbf{P}_{1}$eqv$\mathbf{Q}_{0}\mathbf{Q}_{1}$. The theorem
is proved.

Note, that in the proper Euclidean space, where the concept of equivalence
is single-variant and transitive, the equivalence may be replaced by the
equality, and the relations $\mathbf{P}_{0}\mathbf{Q}_{0}$eqv$\mathbf{P}_{1}%
\mathbf{Q}_{1}$ and $\mathbf{P}_{0}\mathbf{P}_{1}$eqv$\mathbf{Q}_{0}\mathbf{Q%
}_{1}$ may be written respectively in the form $\mathbf{P}_{0}\mathbf{Q}_{0}=%
\mathbf{P}_{1}\mathbf{Q}_{1}$ and $\mathbf{P}_{0}\mathbf{P}_{1}=\mathbf{Q}%
_{0}\mathbf{Q}_{1}$. (Here symbol $"="$ means single-variant equivalence. It
is used in the sense, which it has in the Euclidean geometry.) These
relations are equivalent, and the additional condition $\left\vert \mathbf{P}%
_{0}\mathbf{P}_{1}\right\vert =\left\vert \mathbf{Q}_{0}\mathbf{Q}%
_{1}\right\vert $ is a corollary of any of these relations. Indeed, adding
vector $\mathbf{P}_{1}\mathbf{Q}_{0}$ to both sides of the equality%
\begin{equation}
\mathbf{P}_{0}\mathbf{P}_{1}=\mathbf{Q}_{0}\mathbf{Q}_{1}  \label{a2.11}
\end{equation}%
we obtain 
\begin{equation}
\mathbf{P}_{0}\mathbf{Q}_{0}=\mathbf{P}_{1}\mathbf{Q}_{1}  \label{a2.12}
\end{equation}

Besides, in the proper Euclidean geometry the relation $\mathbf{P}_{1}%
\mathbf{P}_{2}$eqv$\mathbf{Q}_{1}\mathbf{Q}_{2}$ is a corollary of the
relations $\left( \mathbf{P}_{0}\mathbf{P}_{1}\text{eqv}\mathbf{Q}_{0}%
\mathbf{Q}_{1}\right) \wedge \left( \mathbf{P}_{0}\mathbf{P}_{2}\text{eqv}%
\mathbf{Q}_{0}\mathbf{Q}_{2}\right) $. To prove this, we write these
relations in the form of equalities 
\begin{equation}
\mathbf{P}_{0}\mathbf{P}_{1}=\mathbf{Q}_{0}\mathbf{Q}_{1},\qquad \mathbf{P}%
_{0}\mathbf{P}_{2}=\mathbf{Q}_{0}\mathbf{Q}_{2}  \label{a2.14}
\end{equation}%
Subtracting the first relation (\ref{a2.14}) from the second one, we obtain 
\begin{equation}
\mathbf{P}_{1}\mathbf{P}_{2}=\mathbf{Q}_{1}\mathbf{Q}_{2}  \label{a2.15}
\end{equation}%
Thus, in the proper Euclidean geometry among six relations $\mathbf{P}_{0}%
\mathbf{P}_{1}$eqv$\mathbf{Q}_{0}\mathbf{Q}_{1}$, $\mathbf{P}_{0}\mathbf{P}%
_{2}$eqv$\mathbf{Q}_{0}\mathbf{Q}_{2}$, $\mathbf{P}_{1}\mathbf{P}_{2}$eqv$%
\mathbf{Q}_{1}\mathbf{Q}_{2}$ there are only four independent conditions,
whereas in the general case all six conditions are independent, in general.

\section{Examples of equivalent vectors in uniform isotropic multivariant
space-time geometry}

We consider the $\sigma $-space $V_{\mathrm{d}}=\left\{ \sigma ,\mathbb{R}%
^{4}\right\} $ with the world function%
\begin{equation}
\sigma _{\mathrm{d}}=\left\{ 
\begin{array}{c}
\sigma _{\mathrm{M}}+d,\ \ \text{if}\ \ \sigma _{\mathrm{M}}>0 \\ 
\sigma _{\mathrm{M}},\ \ \ \ \text{if\ \ }\sigma _{\mathrm{M}}\leq 0%
\end{array}%
\right. ,\qquad d=\lambda _{0}^{2}=\text{const}>0  \label{a3.1}
\end{equation}%
where $\sigma _{\mathrm{M}}$ is the world function of the $4$-dimensional
space-time of Minkowski. In the inertial coordinate system the world
function $\sigma _{\mathrm{M}}$ has the form 
\begin{equation}
\sigma _{\mathrm{M}}\left( P,P^{\prime }\right) =\sigma _{\mathrm{M}}\left(
x,x^{\prime }\right) =\left( x^{0}-x^{\prime 0}\right) ^{2}-\left( \mathbf{x}%
-\mathbf{x}^{\prime }\right) ^{2}  \label{a3.1a}
\end{equation}%
where coordinates of points $P$ and $P^{\prime }$ are $P=\left\{ x^{0},%
\mathbf{x}\right\} =\left\{ ct,x^{1},x^{2},x^{3}\right\} $, $P^{\prime
}=\left\{ x^{\prime 0},\mathbf{x}^{\prime }\right\} =\left\{ ct^{\prime
},x^{\prime 1},x^{\prime 2},x^{\prime 3}\right\} $ and $c$ is the speed of
the light. The constant $d$ is qualified as a distortion of the distorted
space-time $V_{\mathrm{d}}$, described by the world function $\sigma $. The
constant $\lambda _{0}$ may be considered as an "elementary length"
associated with the distorted space-time $V_{\mathrm{d}}$.

The space-time (\ref{a3.1}) is uniform and isotropic in the sense, that the
world function $\sigma _{\mathrm{d}}$ is invariant with respect to the
simultaneous Poincar\'{e} transformation of both arguments $x$ and $%
x^{\prime }$.

The continual space-time $V_{\mathrm{d}}$ demonstrates evidence of a
discreteness in the sense, that there are no points $x$, $x^{\prime }$,
separated by the timelike interval $\rho =\sqrt{2\sigma \left( x,x^{\prime
}\right) }$, with $\rho \in \left( 0,\lambda _{0}\right) $. It seems rather
unexpected, that the continual space-time may be simultaneously discrete.
Apparently, discreteness of such a kind should be qualified as discreteness
of time.

Let us consider two equivalent timelike vectors $\mathbf{P}_{0}\mathbf{P}%
_{1} $, $\mathbf{Q}_{0}\mathbf{Q}_{1}$. The points $P_{0},P_{1}$, $%
Q_{0},Q_{1}$ have coordinates%
\begin{equation}
P_{0}=\left\{ 0,0,0,0\right\} ,\qquad P_{1}=\left\{ s,0,0,0\right\}
\label{a3.2}
\end{equation}%
\begin{equation}
Q_{0}=\left\{ a,b,0,0\right\} ,\qquad Q_{1}=\left\{ s+a+\alpha _{0},b+\gamma
_{1},\gamma _{2},\gamma _{3}\right\}  \label{a3.3}
\end{equation}%
The time axis is chosen along the vector $\mathbf{P}_{0}\mathbf{P}_{1}$. The
shift vector $\mathbf{P}_{0}\mathbf{Q}_{0}=\left\{ a,b,0,0\right\} $ lies in
the plane of coordinate axes $x^{0}x^{1}$. The length $s$ of the vector $%
\mathbf{P}_{0}\mathbf{P}_{1}$ and parameters $a,b$ of the shift are supposed
to be given. The numbers $\alpha _{0},\beta _{1},\gamma _{2},\gamma _{3}$
are to be determined from the condition that vectors 
\begin{equation}
\mathbf{P}_{0}\mathbf{P}_{1}=\left\{ s,0,0,0\right\} ,\qquad \mathbf{Q}_{0}%
\mathbf{Q}_{1}=\left\{ s+\alpha _{0},\gamma _{1},\gamma _{2},\gamma
_{3}\right\}  \label{a3.3a}
\end{equation}%
are equivalent.

As far the geometry is uniform and isotropic, the relations (\ref{a3.2}), (%
\ref{a3.3}) describe the general case of the points $P_{0},P_{1},Q_{0},Q_{1}$
disposition with timelike vector $\mathbf{P}_{0}\mathbf{P}_{1}$. The
specificity of formulas (\ref{a3.2}), (\ref{a3.3}) is obtained as a result
of proper choice of the coordinate system.

We consider two different cases.

I. All points $P_{0},P_{1},Q_{0},Q_{1}$ are different, 
\begin{equation}
\sigma \left( P,Q\right) =\sigma _{\mathrm{M}}\left( P,Q\right) +\lambda
_{0}^{2},\qquad P\neq Q,\qquad P,Q\in \left\{ P_{0},P_{1},Q_{0},Q_{1}\right\}
\label{a3.4}
\end{equation}

II. $P_{1}=Q_{0},$ all other points $P_{0},P_{1},Q_{1}$ are different, and
all different points are separated by timelike intervals. In this case we
have%
\begin{equation}
\sigma \left( P,Q\right) =\sigma _{\mathrm{M}}\left( P,Q\right) +\lambda
_{0}^{2},\qquad P\neq Q,\qquad P,Q\in \left\{ P_{0},P_{1},Q_{1}\right\}
\label{a3.4a}
\end{equation}%
\begin{equation}
\sigma \left( P_{1},Q_{0}\right) =\sigma _{\mathrm{M}}\left(
P_{1},Q_{0}\right) =0  \label{a3.5}
\end{equation}

In the first case%
\begin{eqnarray}
\left\vert \mathbf{P}_{0}\mathbf{P}_{1}\right\vert ^{2} &=&\left\vert 
\mathbf{P}_{0}\mathbf{P}_{1}\right\vert _{\mathrm{M}}^{2}+2\lambda
_{0}^{2},\qquad \left\vert \mathbf{Q}_{0}\mathbf{Q}_{1}\right\vert
^{2}=\left\vert \mathbf{Q}_{0}\mathbf{Q}_{1}\right\vert _{\mathrm{M}%
}^{2}+2\lambda _{0}^{2}  \label{a3.6} \\
\left\vert \mathbf{P}_{0}\mathbf{Q}_{0}\right\vert ^{2} &=&\left\vert 
\mathbf{P}_{0}\mathbf{Q}_{0}\right\vert _{\mathrm{M}}^{2}+2\lambda
_{0}^{2},\qquad \left\vert \mathbf{P}_{1}\mathbf{Q}_{1}\right\vert
^{2}=\left\vert \mathbf{P}_{1}\mathbf{Q}_{1}\right\vert _{\mathrm{M}%
}^{2}+2\lambda _{0}^{2}  \label{a3.7}
\end{eqnarray}%
Here and in what follows the index "M" means that the quantity is calculated
in the Minkowski space-time. Taking into account definition of the scalar
product (\ref{a1.5}) and relations (\ref{a3.4}), we obtain%
\begin{equation}
\left( \mathbf{P}_{0}\mathbf{P}_{1}.\mathbf{Q}_{0}\mathbf{Q}_{1}\right)
=\left( \mathbf{P}_{0}\mathbf{P}_{1}.\mathbf{Q}_{0}\mathbf{Q}_{1}\right) _{%
\mathrm{M}},\qquad \left( \mathbf{P}_{0}\mathbf{Q}_{0}.\mathbf{P}_{1}\mathbf{%
Q}_{1}\right) =\left( \mathbf{P}_{0}\mathbf{Q}_{0}.\mathbf{P}_{1}\mathbf{Q}%
_{1}\right) _{\mathrm{M}}  \label{a3.8}
\end{equation}

Condition $\mathbf{P}_{0}\mathbf{P}_{1}$eqv$\mathbf{Q}_{0}\mathbf{Q}_{1}$ in
terms of the Minkowski quantities have the form 
\begin{equation}
\left( \mathbf{P}_{0}\mathbf{P}_{1}.\mathbf{Q}_{0}\mathbf{Q}_{1}\right) _{%
\mathrm{M}}=\sqrt{\left( \left\vert \mathbf{P}_{0}\mathbf{P}_{1}\right\vert
_{\mathrm{M}}^{2}+2\lambda _{0}^{2}\right) \left( \left\vert \mathbf{Q}_{0}%
\mathbf{Q}_{1}\right\vert _{\mathrm{M}}^{2}+2\lambda _{0}^{2}\right) }
\label{a3.9}
\end{equation}%
\begin{equation}
\left\vert \mathbf{P}_{0}\mathbf{P}_{1}\right\vert _{\mathrm{M}%
}^{2}=\left\vert \mathbf{Q}_{0}\mathbf{Q}_{1}\right\vert _{\mathrm{M}}^{2}
\label{a3.10}
\end{equation}

Using for $\left( \mathbf{P}_{0}\mathbf{P}_{1}.\mathbf{Q}_{0}\mathbf{Q}%
_{1}\right) _{\mathrm{M}}$ and $\left\vert \mathbf{P}_{0}\mathbf{P}%
_{1}\right\vert _{\mathrm{M}}^{2}$, conventional expression in terms of
coordinates, we obtain instead of (\ref{a3.9}) and (\ref{a3.10}) by means of
(\ref{a3.3a}) 
\begin{equation}
s\left( s+\alpha _{0}\right) =s^{2}+2\lambda _{0}^{2}  \label{a3.11}
\end{equation}%
\begin{equation}
s^{2}=\left( s+\alpha _{0}\right) ^{2}-\gamma _{1}^{2}-\gamma
_{2}^{2}-\gamma _{3}^{2}  \label{a3.12}
\end{equation}%
Solution of these equations gives%
\begin{equation}
\alpha _{0}=\frac{2\lambda _{0}^{2}}{s},\qquad \gamma _{k}=2\lambda _{0}%
\sqrt{1+\frac{\lambda _{0}^{2}}{s^{2}}}\frac{q_{k}}{q},\qquad k=1,2,3
\label{a3.14}
\end{equation}%
where $q_{1},q_{2},q_{3}$ are arbitrary real constants and%
\begin{equation}
q=\sqrt{q_{1}^{2}+q_{2}^{2}+q_{3}^{2}}  \label{a3.15}
\end{equation}%
Thus, coordinates of the point $Q_{1}$%
\begin{equation}
Q_{1}=\left\{ s+a+\frac{2\lambda _{0}^{2}}{s},b+2\lambda _{0}\sqrt{1+\frac{%
\lambda _{0}^{2}}{s^{2}}}\frac{q_{1}}{q},2\lambda _{0}\sqrt{1+\frac{\lambda
_{0}^{2}}{s^{2}}}\frac{q_{2}}{q},2\lambda _{0}\sqrt{1+\frac{\lambda _{0}^{2}%
}{s^{2}}}\frac{q_{3}}{q}\right\}  \label{a3.16}
\end{equation}

In the case, when $\lambda _{0}\ll \left\vert s\right\vert $%
\begin{equation}
Q_{1}\approx \left\{ s+a+\frac{2\lambda _{0}^{2}}{s},b+\frac{2\lambda
_{0}q_{1}}{q},\frac{2\lambda _{0}q_{2}}{q},\frac{2\lambda _{0}q_{3}}{q}%
\right\}  \label{a3.17}
\end{equation}%
The correction, conditioned by distortion, to the coordinate $x^{0}$ is of
the order $\lambda _{0}^{2}$, whereas the correction to other coordinates is
of the order $\lambda _{0}$. Thus, the vector $\mathbf{Q}_{0}\mathbf{Q}_{1}$
is multivariant, and its multivariance is described by two arbitrary
parameters.

The second case is more interesting from physical viewpoint, because it may
be considered as a description of the temporal evolution of the geometrical
object, described by two characteristic points $P_{0}$ and $P_{1}$,
separated by the timelike interval. In the second case, when $P_{1}=Q_{0}$,
we have $a=s$, $b=0$, i.e.%
\begin{equation}
P_{0}=\left\{ 0,0,0,0\right\} ,\qquad Q_{0}=P_{1}=\left\{ s,0,0,0\right\}
,\qquad Q_{1}=\left\{ 2s+\alpha _{0},\gamma _{1},\gamma _{2},\gamma
_{3}\right\}  \label{a3.18}
\end{equation}%
The vectors $\mathbf{P}_{0}\mathbf{P}_{1}$ and $\mathbf{Q}_{0}\mathbf{Q}_{1}$
have the same form (\ref{a3.3a}), however in this case the relation between $%
\left( \mathbf{P}_{0}\mathbf{P}_{1}.\mathbf{Q}_{0}\mathbf{Q}_{1}\right) $
and $\left( \mathbf{P}_{0}\mathbf{P}_{1}.\mathbf{P}_{1}\mathbf{Q}_{1}\right)
_{\mathrm{E}}$ has the form 
\begin{equation}
\left( \mathbf{P}_{0}\mathbf{P}_{1}.\mathbf{Q}_{0}\mathbf{Q}_{1}\right)
=\left( \mathbf{P}_{0}\mathbf{P}_{1}.\mathbf{P}_{1}\mathbf{Q}_{1}\right)
=\left( \mathbf{P}_{0}\mathbf{P}_{1}.\mathbf{P}_{1}\mathbf{Q}_{1}\right) _{%
\mathrm{M}}-\lambda _{0}^{2},  \label{a3.19}
\end{equation}%
which distinguishes from the relation (\ref{a3.8}).

Condition $\mathbf{P}_{0}\mathbf{P}_{1}$eqv$\mathbf{Q}_{0}\mathbf{Q}_{1}$ in
terms of Minkowskian quantities have the form%
\begin{eqnarray}
\left\vert \mathbf{P}_{0}\mathbf{Q}_{0}\right\vert ^{2} &=&\left\vert 
\mathbf{P}_{0}\mathbf{P}_{1}\right\vert ^{2}=\left\vert \mathbf{P}_{0}%
\mathbf{P}_{1}\right\vert _{\mathrm{M}}^{2}+2\lambda _{0}^{2},\qquad
\left\vert \mathbf{P}_{1}\mathbf{Q}_{0}\right\vert ^{2}=\left\vert \mathbf{P}%
_{1}\mathbf{Q}_{0}\right\vert _{\mathrm{M}}^{2}=0  \label{a3.20} \\
\left\vert \mathbf{P}_{1}\mathbf{Q}_{1}\right\vert ^{2} &=&\left\vert 
\mathbf{Q}_{0}\mathbf{Q}_{1}\right\vert ^{2}=\left\vert \mathbf{Q}_{0}%
\mathbf{Q}_{1}\right\vert _{\mathrm{M}}^{2}+2\lambda _{0}^{2}  \label{a3.21}
\end{eqnarray}%
\begin{equation}
\left( \mathbf{P}_{0}\mathbf{P}_{1}.\mathbf{Q}_{0}\mathbf{Q}_{1}\right) _{%
\mathrm{M}}-\lambda _{0}^{2}=\sqrt{\left( \left\vert \mathbf{P}_{0}\mathbf{P}%
_{1}\right\vert _{\mathrm{M}}^{2}+2\lambda _{0}^{2}\right) \left( \left\vert 
\mathbf{Q}_{0}\mathbf{Q}_{1}\right\vert _{\mathrm{M}}^{2}+2\lambda
_{0}^{2}\right) }  \label{a3.29}
\end{equation}%
\begin{equation}
\left\vert \mathbf{P}_{0}\mathbf{P}_{1}\right\vert _{\mathrm{M}%
}^{2}=\left\vert \mathbf{Q}_{0}\mathbf{Q}_{1}\right\vert _{\mathrm{M}}^{2}
\label{a3.30}
\end{equation}%
They take the form%
\begin{equation}
s\left( s+\alpha _{0}\right) -\lambda _{0}^{2}=s^{2}+2\lambda _{0}^{2}
\label{a3.31}
\end{equation}%
\begin{equation}
s^{2}=\left( s+\alpha _{0}\right) ^{2}-\gamma _{1}^{2}-\gamma
_{2}^{2}-\gamma _{3}^{2}  \label{a3.32}
\end{equation}%
and distinguish from the relations (\ref{a3.11}), (\ref{a3.12}) by
replacement of $\lambda _{0}$ by $\sqrt{3/2}\lambda _{0}$.

Using the change $a\rightarrow s$, $b\rightarrow 0$, $\lambda _{0}$ $%
\rightarrow \sqrt{3/2}\lambda _{0}$ in relations (\ref{a3.16}), (\ref{a3.17}%
), we obtain for the vector $\mathbf{Q}_{0}\mathbf{Q}_{1}$%
\begin{equation}
\mathbf{Q}_{0}\mathbf{Q}_{1}=\mathbf{P}_{1}\mathbf{Q}_{1}=\left\{ s+\frac{%
3\lambda _{0}^{2}}{s},\lambda _{0}\kappa \frac{q_{1}}{q},\lambda _{0}\kappa 
\frac{q_{2}}{q},\lambda _{0}\kappa \frac{q_{3}}{q}\right\}  \label{a3.33}
\end{equation}%
where $q_{1},q_{2},q_{3}$ are arbitrary quantities%
\begin{equation}
\kappa =\sqrt{6\left( 1+\frac{3\lambda _{0}^{2}}{2s^{2}}\right) }
\label{a3.33a}
\end{equation}%
and $q$ is determined by the relation (\ref{a3.15}).

In the case, when $\lambda _{0}\ll \left\vert s\right\vert $%
\begin{equation}
\mathbf{Q}_{0}\mathbf{Q}_{1}=\mathbf{P}_{1}\mathbf{Q}_{1}\approx \left\{ s+%
\frac{3\lambda _{0}^{2}}{s},\frac{\sqrt{6}\lambda _{0}q_{1}}{q},\frac{\sqrt{6%
}\lambda _{0}q_{2}}{q},\frac{\sqrt{6}\lambda _{0}q_{3}}{q}\right\}
\label{a3.34}
\end{equation}

Another limit case, when $s\ll \lambda _{0}$. We set $s=\beta \lambda _{0}$, 
$\beta \ll 1$ and obtain%
\begin{equation}
\mathbf{Q}_{0}\mathbf{Q}_{1}=\mathbf{P}_{1}\mathbf{Q}_{1}=\left\{ \lambda
_{0}\left( \beta +\frac{3}{\beta }\right) ,\lambda _{0}\kappa _{1}\frac{q_{1}%
}{q},\lambda _{0}\kappa _{1}\frac{q_{2}}{q},\lambda \kappa _{1}\frac{q_{3}}{q%
}\right\}  \label{a3.34a}
\end{equation}%
where%
\begin{equation}
\kappa _{1}=\sqrt{6\left( 1+\frac{3}{2\beta ^{2}}\right) }  \label{a3.35}
\end{equation}%
In this case all coefficients $\alpha _{0},\gamma _{1},\gamma _{2},\gamma
_{3}$ are of the same order $3\lambda _{0}/\beta $, and multivariance of the
vector $\mathbf{Q}_{0}\mathbf{Q}_{1}=\mathbf{P}_{1}\mathbf{Q}_{1}$ is very
large. We have for the $\left\vert \mathbf{Q}_{0}\mathbf{Q}_{1}\right\vert
^{2}$%
\begin{equation}
\left\vert \mathbf{Q}_{0}\mathbf{Q}_{1}\right\vert ^{2}=\left\vert \mathbf{Q}%
_{0}\mathbf{Q}_{1}\right\vert _{\mathrm{M}}^{2}+2\lambda _{0}^{2}=\left(
2+\beta ^{2}\right) \lambda _{0}^{2}  \label{a3.36}
\end{equation}

\section{Equivalence of two null vectors}

Let us consider two equivalent null vectors $\mathbf{P}_{0}\mathbf{P}_{1}$
and $\mathbf{Q}_{0}\mathbf{Q}_{1}$ in the space-time (\ref{a3.1}) $%
\left\vert \mathbf{P}_{0}\mathbf{P}_{1}\right\vert =0,$ $\left\vert \mathbf{Q%
}_{0}\mathbf{Q}_{1}\right\vert =0$. We consider the case, when the points $%
P_{1}$ and $Q_{0}$ coincide $P_{1}=Q_{0}$. Then 
\begin{equation}
P_{0}=\left\{ 0,0,0,0\right\} ,\qquad P_{1}=\left\{ s,s,0,0\right\} ,
\label{b5.10}
\end{equation}%
\begin{subequations}
\begin{equation}
Q_{0}=\left\{ s,s,0,0\right\} ,\qquad Q_{1}=\left\{ 2s+\alpha _{0},2s+\gamma
_{1},\gamma _{2},\gamma _{3}\right\} ,  \label{b5.11}
\end{equation}%
\end{subequations}
\begin{equation}
\mathbf{P}_{0}\mathbf{P}_{1}=\left\{ s,s,0,0\right\} ,\qquad \mathbf{P}_{1}%
\mathbf{Q}_{1}=\mathbf{Q}_{0}\mathbf{Q}_{1}=\left\{ s+\alpha _{0},s+\gamma
_{1},\gamma _{2},\gamma _{3}\right\}   \label{b5.12}
\end{equation}%
\begin{equation}
\mathbf{P}_{0}\mathbf{Q}_{1}=\left\{ 2s+\alpha _{0},2s+\gamma _{1},\gamma
_{2},\gamma _{3}\right\}   \label{b5.12a}
\end{equation}

In this case we obtain 
\begin{equation}
\left( \mathbf{P}_{0}\mathbf{P}_{1}.\mathbf{Q}_{0}\mathbf{Q}_{1}\right)
=\sigma \left( P_{0},Q_{1}\right) +\sigma \left( P_{1},Q_{0}\right) -\sigma
\left( P_{0},Q_{0}\right) -\sigma \left( P_{1},Q_{1}\right) =\sigma \left(
P_{0},Q_{1}\right)   \label{c5.14}
\end{equation}%
The equivalence conditions have the form%
\begin{equation}
\left\vert \mathbf{P}_{0}\mathbf{P}_{1}\right\vert ^{2}=\left\vert \mathbf{P}%
_{0}\mathbf{P}_{1}\right\vert _{\mathrm{M}}^{2}=0,\qquad \left\vert \mathbf{P%
}_{1}\mathbf{Q}_{2}\right\vert ^{2}=\left\vert \mathbf{P}_{1}\mathbf{Q}%
_{2}\right\vert _{\mathrm{M}}^{2}=0  \label{a7.22}
\end{equation}%
\begin{equation}
\left( \mathbf{P}_{0}\mathbf{P}_{1}.\mathbf{P}_{1}\mathbf{Q}_{2}\right)
=\sigma \left( P_{0},Q_{2}\right) -\sigma \left( P_{1},Q_{2}\right) -\sigma
\left( P_{0},P_{1}\right) =\sigma \left( P_{0},Q_{2}\right) =0  \label{a7.23}
\end{equation}

In terms of coordinates we have%
\begin{equation}
\sigma \left( P_{0},Q_{2}\right) =\sigma _{\mathrm{M}}\left(
P_{0},Q_{2}\right) =\left( 2s+\alpha _{0}\right) ^{2}-\left( 2s+\gamma
_{1}\right) ^{2}-\gamma _{2}^{2}-\gamma _{3}^{2}=0  \label{a7.24}
\end{equation}%
\begin{equation}
\left( s+\alpha _{0}\right) ^{2}-\left( s+\gamma _{1}\right) ^{2}-\gamma
_{2}^{2}-\gamma _{3}^{3}=  \label{a7.25}
\end{equation}%
Solution of equations (\ref{a7.24}), (\ref{a7.25}) has the form%
\begin{equation*}
\alpha _{0}=\gamma _{1},\qquad \gamma _{2}=\gamma _{3}=0
\end{equation*}%
\begin{equation}
\mathbf{P}_{1}\mathbf{Q}_{2}=\left\{ s+\alpha _{0},s+\alpha _{0},0,0\right\}
,\qquad \mathbf{P}_{0}\mathbf{Q}_{2}=\left\{ 2s+\alpha _{0},2s+\alpha
_{0},0,0\right\}   \label{a7.26}
\end{equation}%
where $\alpha _{0}$ is an arbitrary real quantity.

\section{Construction of geometrical objects in T-geometry}

Geometrical object $\mathcal{O\subset }\Omega $ is a subset of points in the
point space $\Omega $. In the T-geometry the geometric object $\mathcal{O}$
is described by means of the skeleton-envelope method \cite{R02}. It means
that any geometric object $\mathcal{O}$ is considered to be a set of
intersections and joins of elementary geometric objects (EGO).

The finite set $\mathcal{P}^{n}\equiv \left\{ P_{0},P_{1},...,P_{n}\right\}
\subset \Omega $ of parameters of the envelope function $f_{\mathcal{P}^{n}}$
is the skeleton of elementary geometric object (EGO) $\mathcal{E}\subset
\Omega $. The set $\mathcal{E}\subset \Omega $ of points forming EGO is
called the envelope of its skeleton $\mathcal{P}^{n}$. In the continuous
generalized geometry the envelope $\mathcal{E}$ is usually a continual set
of points. The envelope function $f_{\mathcal{P}^{n}}$%
\begin{equation}
f_{\mathcal{P}^{n}}:\qquad \Omega \rightarrow \mathbb{R},  \label{h2.1}
\end{equation}%
determining EGO, is a function of the running point $R\in \Omega $ and of
parameters $\mathcal{P}^{n}\subset \Omega $. The envelope function $f_{%
\mathcal{P}^{n}}$ is supposed to be an algebraic function of $s$ arguments $%
w=\left\{ w_{1},w_{2},...w_{s}\right\} $, $s=(n+2)(n+1)/2$. Each of
arguments $w_{k}=\sigma \left( Q_{k},L_{k}\right) $ is the world function $%
\sigma $ of two arguments $Q_{k},L_{k}\in \left\{ R,\mathcal{P}^{n}\right\} $%
, either belonging to skeleton $\mathcal{P}^{n}$, or coinciding with the
running point $R$. Thus, any elementary geometric object $\mathcal{E}$ is
determined by its skeleton $\mathcal{P}^{n}$ and its envelope function $f_{%
\mathcal{P}^{n}}$ as the set of zeros of the envelope function 
\begin{equation}
\mathcal{E}=\left\{ R|f_{\mathcal{P}^{n}}\left( R\right) =0\right\}
\label{h2.2}
\end{equation}

Characteristic points of the EGO are the skeleton points $\mathcal{P}%
^{n}\equiv \left\{ P_{0},P_{1},...,P_{n}\right\} $. The simplest example of
EGO is the segment $\mathcal{T}_{\left[ P_{0}P_{1}\right] }$ of the straight
line (\ref{a1.17}) between the points $P_{0}$ and $P_{1}$, which is defined
by the relation 
\begin{eqnarray}
\mathcal{T}_{\left[ P_{0}P_{1}\right] } &=&\left\{ R|f_{P_{0}P_{1}}\left(
R\right) =0\right\} ,  \notag \\
f_{P_{0}P_{1}}\left( R\right) &=&\sqrt{2\sigma \left( P_{0},R\right) }+\sqrt{%
2\sigma \left( R,P_{1}\right) }-\sqrt{2\sigma \left( P_{0},P_{1}\right) }
\label{a5.1}
\end{eqnarray}

Another example is the sphere $\mathcal{S}_{OQ}$, where $O$ is the center of
the sphere and $Q$ is some point on the surface of the sphere. The sphere $%
\mathcal{S}_{P_{0}P_{1}}$ is described by the relation 
\begin{equation}
\mathcal{S}_{OQ}=\left\{ R|g_{OQ}\left( R\right) =0\right\} ,\qquad
g_{OQ}\left( R\right) =\sqrt{2\sigma \left( O,R\right) }-\sqrt{2\sigma
\left( O,Q\right) }  \label{a5.2}
\end{equation}%
Here points $O,Q$ form the skeleton of the sphere, whereas the function $%
g_{OQ}$ is the envelope function.

The third example is the cylinder $\mathcal{C}(P_{0},P_{1},Q)$ with the
points $P_{0},P_{1}$ on the cylinder axis and the point $Q$ on its surface.
The cylinder $\mathcal{C}(P_{0},P_{1},Q)$ is determined by the relation 
\begin{eqnarray}
\mathcal{C}(P_{0},P_{1},Q) &=&\left\{ R|f_{P_{0}P_{1}Q}\left( R\right)
=0\right\} ,  \label{g3.1} \\
f_{P_{0}P_{1}Q}\left( R\right) &=&F_{2}\left( P_{0},P_{1},Q\right)
-F_{2}\left( P_{0},P_{1},R\right)  \notag
\end{eqnarray}%
\begin{equation}
F_{2}\left( P_{0},P_{1},Q\right) =\left\vert 
\begin{array}{cc}
\left( \mathbf{P}_{0}\mathbf{P}_{1}.\mathbf{P}_{0}\mathbf{P}_{1}\right) & 
\left( \mathbf{P}_{0}\mathbf{P}_{1}.\mathbf{P}_{0}\mathbf{Q}\right) \\ 
\left( \mathbf{P}_{0}\mathbf{Q}.\mathbf{P}_{0}\mathbf{P}_{1}\right) & \left( 
\mathbf{P}_{0}\mathbf{Q}.\mathbf{P}_{0}\mathbf{Q}\right)%
\end{array}%
\right\vert  \label{g3.2}
\end{equation}%
Here $\sqrt{F_{2}\left( P_{0},P_{1},Q\right) }$ is the area of the
parallelogram, constructed on the vectors $\mathbf{P}_{0}\mathbf{P}_{1}$ and 
$\mathbf{P}_{0}\mathbf{Q}$\textbf{\ }and $\frac{1}{2}\sqrt{F_{2}\left(
P_{0},P_{1},Q\right) }\ $ is the area of triangle with vertices at the
points $P_{0},P_{1},Q$. The equality $F_{2}\left( P_{0},P_{1},Q\right)
=F_{2}\left( P_{0},P_{1},R\right) $ means that the distance between the
point $Q$ and the axis, determined by the vector $\mathbf{P}_{0}\mathbf{P}%
_{1}$, is equal to the distance between $R$ and the axis. Here the points $%
P_{0},P_{1},Q$ form the skeleton of the cylinder, whereas the function $%
f_{P_{0}P_{1}Q}$ is the envelope function.

\textit{Definition.} Two EGOs $\mathcal{E}\left( \mathcal{P}^{n}\right) $
and $\mathcal{E}\left( \mathcal{Q}^{n}\right) $ are equivalent, if their
skeletons are equivalent and their envelope functions $f_{\mathcal{P}^{n}}$
and $g_{\mathcal{Q}^{n}}$ are equal. Equivalence of two skeletons $\mathcal{P%
}^{n}\equiv \left\{ P_{0},P_{1},...,P_{n}\right\} \subset \Omega $ and $%
\mathcal{Q}^{n}\equiv \left\{ Q_{0},Q_{1},...,Q_{n}\right\} \subset \Omega $
means that 
\begin{equation}
\mathbf{P}_{i}\mathbf{P}_{k}\text{eqv}\mathbf{Q}_{i}\mathbf{Q}_{k},\qquad
i,k=0,1,...n,\quad i<k  \label{a5.4}
\end{equation}%
Equivalence of the envelope functions $f_{\mathcal{P}^{n}}$ and $g_{\mathcal{%
Q}^{n}}$ means that 
\begin{equation}
f_{\mathcal{P}^{n}}\left( R\right) =\Phi \left( g_{\mathcal{P}^{n}}\left(
R\right) \right) ,\qquad \forall R\in \Omega  \label{a5.5}
\end{equation}%
where $\Phi $ is an arbitrary function, having the property%
\begin{equation}
\Phi :\mathbb{R}\rightarrow \mathbb{R},\qquad \Phi \left( 0\right) =0
\label{a5.5a}
\end{equation}

Equivalence of shapes of two EGOs $\mathcal{E}\left( \mathcal{P}^{n}\right) $
and $\mathcal{E}\left( \mathcal{Q}^{n}\right) $ is determined by equivalence
of shapes of their skeletons $\mathcal{P}^{n}$ and $\mathcal{Q}^{n}$, which
is described by the relations

\begin{equation}
\left\vert \mathbf{P}_{i}\mathbf{P}_{k}\right\vert =\left\vert \mathbf{Q}_{i}%
\mathbf{Q}_{k}\right\vert ,\qquad i,k=0,1,...n,\quad i<k  \label{a5.6}
\end{equation}%
and equivalence of their envelope functions $f_{\mathcal{P}^{n}}$ and $g_{%
\mathcal{Q}^{n}}$ (\ref{a5.5}).

Equivalence of orientations of skeletons $\mathcal{P}^{n}$ and $\mathcal{Q}%
^{n}$ in the point space $\Omega $ is described by the relations%
\begin{equation}
\mathbf{P}_{i}\mathbf{P}_{k}\upuparrows \mathbf{Q}_{i}\mathbf{Q}_{k},\qquad
i,k=0,1,...n,\quad i<k  \label{a5.7}
\end{equation}

Equivalence of shapes and orientations of skeletons is equivalence of
skeletons, described by the relations (\ref{a5.4}).

\textit{Definition}. The elementary geometric object $\mathcal{E}\left( 
\mathcal{P}^{n}\right) $ exists, if at any time moment there is an
elementary geometric object $\mathcal{E}^{\prime }\left( \mathcal{P}^{\prime
n}\right) $, which is equivalent to EGO $\mathcal{E}\left( \mathcal{P}%
^{n}\right) $. We suppose, the skeleton $\mathcal{P}^{n}=\left\{
P_{0},P_{1},...P_{n}\right\} $ contains points separated by the timelike
interval. We suppose that the vector $\mathbf{P}_{0}\mathbf{P}_{1}$ is
timelike ($\left\vert \mathbf{P}_{0}\mathbf{P}_{1}\right\vert ^{2}>0$). We
assume, that the elementary geometrical object $\mathcal{E}\left( \mathcal{P}%
^{n}\right) =\mathcal{E}\left( P_{0},P_{1},...,P_{n}\right) $ is placed at
the point $P_{0}$. The same EGO placed at the point $P_{1}$ has the form $%
\mathcal{E}\left( P_{0}^{\prime },P_{1}^{\prime },P_{2}^{\prime
}...,P_{n}^{\prime }\right) $ with $P_{0}^{\prime }=P_{1}$. The points $%
P_{0} $ and $P_{0}^{\prime }=P_{1}$ are separated by the timelike interval,
and we may consider the EGO $\mathcal{E}\left( P_{0}^{\prime },P_{1}^{\prime
},P_{2}^{\prime }...,P_{n}^{\prime }\right) $ as a result of temporal
evolution of the EGO $\mathcal{E}\left( P_{0},P_{1},...,P_{n}\right) $,
provided these objects are equivalent, i.e. 
\begin{equation}
\mathbf{P}_{i}\mathbf{P}_{k}\text{eqv}\mathbf{P}_{i}^{\prime }\mathbf{P}%
_{k}^{\prime },\qquad i,k=0,1,...n,\quad i<k  \label{a6.1}
\end{equation}

Thus, if EGOs $\mathcal{E}\left( P_{0}^{\left( 0\right) },P_{1}^{\left(
0\right) },...,P_{n}^{\left( 0\right) }\right) $, $\mathcal{E}\left(
P_{0}^{\left( 1\right) },P_{1}^{\left( 1\right) },...,P_{n}^{\left( 1\right)
}\right) $, $\mathcal{E}\left( P_{0}^{\left( 2\right) },P_{1}^{\left(
2\right) },...,P_{n}^{\left( 2\right) }\right) $,... $\mathcal{E}\left(
P_{0}^{\left( k\right) },P_{1}^{\left( k\right) },...,P_{n}^{\left( k\right)
}\right) $,... are equivalent in pairs 
\begin{equation}
\mathcal{E}\left( P_{0}^{\left( k-1\right) },P_{1}^{\left( k-1\right)
},...,P_{n}^{\left( k-1\right) }\right) \text{eqv}\mathcal{E}\left(
P_{0}^{\left( k\right) },P_{1}^{\left( k\right) },...,P_{n}^{\left( k\right)
}\right) ,\qquad k=1,2,...  \label{a6.2}
\end{equation}%
and 
\begin{eqnarray}
P_{1}^{\left( 0\right) } &=&P_{0}^{\left( 1\right) },\quad P_{1}^{\left(
1\right) }=P_{0}^{\left( 2\right) },\quad P_{1}^{\left( 2\right)
}=P_{0}^{\left( 3\right) },...P_{1}^{\left( k\right) }=P_{0}^{\left(
k+1\right) },...  \label{a6.3} \\
\left\vert \mathbf{P}_{0}^{\left( k\right) }\mathbf{P}_{1}^{\left( k\right)
}\right\vert ^{2} &>&0,\qquad k=1,2,...  \label{a6.4}
\end{eqnarray}%
one may consider existence of the set of elementary geometrical objects
(EGOs) $\mathcal{E}\left( P_{0}^{\left( 0\right) },P_{1}^{\left( 0\right)
},...,P_{n}^{\left( 0\right) }\right) $, $\mathcal{E}\left( P_{0}^{\left(
1\right) },P_{1}^{\left( 1\right) },...,P_{n}^{\left( 1\right) }\right) $, $%
\mathcal{E}\left( P_{0}^{\left( 2\right) },P_{1}^{\left( 2\right)
},...,P_{n}^{\left( 2\right) }\right) $,...\newline
$\mathcal{E}\left( P_{0}^{\left( k\right) },P_{1}^{\left( k\right)
},...,P_{n}^{\left( k\right) }\right) $, with properties (\ref{a6.2}) - (\ref%
{a6.4}) as a temporal evolution of EGO $\mathcal{E}\left( P_{0}^{\left(
0\right) },P_{1}^{\left( 0\right) },...,P_{n}^{\left( 0\right) }\right) $.

Thus, the space-time geometry determines a possibility of existence of the
geometric object $\mathcal{E}\left( \mathcal{P}^{n}\right) $ and its
temporal evolution. Some objects may exist, other ones do not exist. This
fact depends on possibility of fulfilment of the relation (\ref{a6.2}). For
some geometrical objects the temporal evolution may be multivariant, for
other ones it is single-variant. It is possible such geometrical objects,
for which there is no equivalent geometrical objects, and there is no
temporal evolution.

\section{Temporal evolution of timelike segment of the straight}

Timelike segment (\ref{a5.1}) of the straight line (\ref{a1.17}) is an
elementary geometrical object $\mathcal{T}_{\left[ P_{0}P_{1}\right] }$,
described by the skeleton, consisting of two points $P_{0},P_{1}$. Temporal
evolution of this segment is described by the broken tube $\mathcal{T}_{%
\mathrm{br}}$. 
\begin{equation}
\mathcal{T}_{\mathrm{br}}=\dbigcup\limits_{i}\mathcal{T}_{\left[ P_{i}P_{i+1}%
\right] },\qquad \mathbf{P}_{i-1}\mathbf{P}_{i}\text{eqv}\mathbf{P}_{i}%
\mathbf{P}_{i+1},\qquad \left\vert \mathbf{P}_{i}\mathbf{P}_{i+1}\right\vert
^{2}=\mu ^{2},\qquad i=0,\pm 1,\pm 2,...  \label{b6.1}
\end{equation}%
The shape of $\mathcal{T}_{\mathrm{br}}$ in the space-time (\ref{a3.1}) is
multivariant. It looks as a chain, consisting of similar links $\mathcal{T}_{%
\left[ P_{i}P_{i+1}\right] }$. In the Minkowski space-time the shape of $%
\mathcal{T}_{\mathrm{br}}$ is single-variant, and the chain of links $%
\mathcal{T}_{\left[ P_{i}P_{i+1}\right] }$ degenerates into the timelike
straight line.

In the inertial coordinate system $\left\{ ct,x^{1},x^{2},x^{3}\right\} $,
where the points $P_{0},P_{1}$ have coordinates%
\begin{equation}
P_{0}=\left\{ 0,0,0,0\right\} ,\qquad P_{1}=\left\{ \sqrt{\mu ^{2}-2\lambda
_{0}^{2}},0,0,0\right\} ,  \label{b6.2}
\end{equation}%
and the vector $\mathbf{P}_{0}\mathbf{P}_{1}$ has the length%
\begin{equation}
\left\vert \mathbf{P}_{0}\mathbf{P}_{1}\right\vert =\sqrt{\left\vert \mathbf{%
P}_{0}\mathbf{P}_{1}\right\vert _{\mathrm{M}}^{2}+2\lambda _{0}^{2}}=\mu ,
\label{b6.4}
\end{equation}%
the surface $\mathcal{T}_{\left( P_{0}P_{1}\right) }=\mathcal{T}_{\left[
P_{0}P_{1}\right] }\backslash \left\{ P_{0},P_{1}\right\} $ is described by
the following equation 
\begin{equation}
r^{2}=\mathbf{x}^{2}=\frac{2\lambda _{0}^{2}\left( ct-\frac{1}{2}\sqrt{\mu
^{2}-2\lambda _{0}^{2}}\right) ^{2}}{\mu ^{2}}+\frac{3}{2}\lambda
_{0}^{2},\qquad 0<ct<\sqrt{\mu ^{2}-2\lambda _{0}^{2}}  \label{b6.3}
\end{equation}%
The skeleton points $P_{0},P_{1}$ do not belong to the surface (\ref{b6.3}),
but the belong to $\mathcal{T}_{\left[ P_{0}P_{1}\right] }$. This surface is
a tube with minimal radius $r_{\mathrm{\min }}=\sqrt{3/2}\lambda _{0}$ and
maximal radius $r_{\mathrm{\max }}=\sqrt{2\lambda _{0}^{2}-\frac{\lambda
_{0}^{4}}{\mu ^{2}}}$, $\mu >\sqrt{2}\lambda _{0}$. In the limit $\lambda
_{0}\rightarrow 0$ the maximal and minimal radii tend to zero, and the tube (%
\ref{b6.3}) degenerates into the straight line interval.

\textit{Remark. }The form of the envelope function is of no importance for
the temporal evolution of the geometrical object, described by two skeleton
points $P_{0},P_{1}$, In particular, one may consider the sphere $\mathcal{S}%
_{P_{0}P_{1}}$, defined by the relation (\ref{a5.2}), instead of EGO $%
\mathcal{T}_{\left[ P_{0}P_{1}\right] }$. In this case we have instead of
the broken tube (\ref{b6.1}) 
\begin{equation}
\mathcal{T}_{\mathrm{br}}=\dbigcup\limits_{i}\mathcal{S}_{P_{i}P_{i+1}},%
\qquad \mathbf{P}_{i-1}\mathbf{P}_{i}\text{eqv}\mathbf{P}_{i}\mathbf{P}%
_{i+1},\qquad \left\vert \mathbf{P}_{i}\mathbf{P}_{i+1}\right\vert ^{2}=\mu
^{2},\qquad i=0,\pm 1,\pm 2,...  \label{b6.5}
\end{equation}%
Here instead of the tube segment (\ref{b6.3}), we have the hyperbola%
\begin{equation}
r^{2}=\mathbf{x}^{2}=c^{2}t^{2}-\mu ^{2}+2\lambda _{0}^{2}  \label{b6.6}
\end{equation}%
Defect of the presentation (\ref{b6.5}) is determined by the fact that one
of the skeleton points ($P_{0}$) does not belong to the sphere envelope $%
\mathcal{S}_{P_{0}P_{1}}$. In the limit $\lambda _{0}\rightarrow 0$ the set
of hyperbolas (\ref{b6.6}) is not associated with the particle world line.

The mutual location of two adjacent links $\mathcal{T}_{\left[ P_{0}P_{1}%
\right] }$ and $\mathcal{T}_{\left[ P_{1}P_{2}\right] }$ is described by the
angle between the vectors $\mathbf{P}_{0}\mathbf{P}_{1}$ and $\mathbf{P}_{1}%
\mathbf{P}_{2}$. As far as these vectors are equivalent and, hence, are in
parallel this angle $\theta =0$, because 
\begin{equation}
\cosh \theta =\frac{\left( \mathbf{P}_{0}\mathbf{P}_{1}.\mathbf{P}_{1}%
\mathbf{P}_{2}\right) }{\left\vert \mathbf{P}_{0}\mathbf{P}_{1}\right\vert
\cdot \left\vert \mathbf{P}_{1}\mathbf{P}_{2}\right\vert }=1  \label{b6.7}
\end{equation}%
However, as far as 
\begin{equation}
\left\vert \mathbf{P}_{0}\mathbf{P}_{1}\right\vert _{\mathrm{M}%
}^{2}=\left\vert \mathbf{P}_{1}\mathbf{P}_{2}\right\vert _{\mathrm{M}%
}^{2}=\left\vert \mathbf{P}_{0}\mathbf{P}_{1}\right\vert ^{2}-2\lambda
_{0}^{2},\qquad \left( \mathbf{P}_{0}\mathbf{P}_{1}.\mathbf{P}_{1}\mathbf{P}%
_{2}\right) _{\mathrm{M}}=\left( \mathbf{P}_{0}\mathbf{P}_{1}.\mathbf{P}_{1}%
\mathbf{P}_{2}\right) -\lambda _{0}^{2}  \label{b6.8}
\end{equation}%
the angle $\theta _{\mathrm{M}}$ between vectors $\mathbf{P}_{0}\mathbf{P}%
_{1}$ and $\mathbf{P}_{1}\mathbf{P}_{2}$, measured on the Minkowski manifold
is defined by the relation 
\begin{equation}
\cosh \theta _{\mathrm{M}}=\frac{\left( \mathbf{P}_{0}\mathbf{P}_{1}.\mathbf{%
P}_{1}\mathbf{P}_{2}\right) _{\mathrm{M}}}{\left\vert \mathbf{P}_{0}\mathbf{P%
}_{1}\right\vert _{\mathrm{M}}\cdot \left\vert \mathbf{P}_{1}\mathbf{P}%
_{2}\right\vert _{\mathrm{M}}}=\frac{\left( \mathbf{P}_{0}\mathbf{P}_{1}.%
\mathbf{P}_{1}\mathbf{P}_{2}\right) -\lambda _{0}^{2}}{\left\vert \mathbf{P}%
_{0}\mathbf{P}_{1}\right\vert ^{2}-2\lambda _{0}^{2}}=\frac{\mu ^{2}-\lambda
_{0}^{2}}{\mu ^{2}-2\lambda _{0}^{2}}>1  \label{b6.9}
\end{equation}

In the case, when $\mu \gg \lambda _{0}$, it follows from (\ref{b6.9}), that 
\begin{equation}
\theta _{\mathrm{M}}=\sqrt{2}\frac{\lambda _{0}}{\mu }  \label{b6.10}
\end{equation}%
Thus, at fixed vector $\mathbf{P}_{0}\mathbf{P}_{1}$ the point $P_{2}$,
determining the adjacent vector $\mathbf{P}_{1}\mathbf{P}_{2}$, lies on the
surface of the cone with angle $\theta _{\mathrm{M}}$ at the vertex $P_{1}$.
Thus, if $\lambda _{0}\neq 0$, the broken tube (\ref{b6.1}) is multivariant.
If $\lambda _{0}\rightarrow 0$, then according to (\ref{b6.10}) $\theta _{%
\mathrm{M}}\rightarrow 0$, and the cone degenerates into the straight line.
The broken tube (\ref{b6.1}) degenerates into the straight line, which
associates with the world line of a particle. In the Minkowski space-time
the world line is a geometric characteristic of a particle. However,
besides, the particle has dynamic characteristics: the momentum $\mathbf{p}$
and the mass, which cannot be determined geometrically in the Minkowski
space-time. The momentum vector is tangent to the world line, and its
direction may be determined by the shape of the world line. However, its
length (the particle mass) cannot be determined geometrically. In the
Minkowski space-time the mass is a non-geometric characteristic, associated
with the particle world line.

In the distorted space-time $V_{\mathrm{d}}$, described by the world
function (\ref{a3.1}), the particle mass is defined geometrically as the
length $\left\vert \mathbf{P}_{i}\mathbf{P}_{i+1}\right\vert =\mu $ of the
link of the broken tube (\ref{b6.1}). Indeed, in $V_{\mathrm{d}}$ the link $%
\mathcal{T}_{\left[ P_{i}P_{i+1}\right] }$ is a tube segment of the radius $%
r $. The ends $P_{i}$ and $P_{i+1}$ of this segment can be determined
geometrically. If we know the ends $P_{i}$, $P_{i+1}$ of $\mathcal{T}_{\left[
P_{i}P_{i+1}\right] }$, we can determine the length $\left\vert \mathbf{P}%
_{i}\mathbf{P}_{i+1}\right\vert =\mu $. But in this case the mass $\mu $ is
determined in the units of length as the distance between the skeleton
points $P_{i},P_{i+1}$ of the segment $\mathcal{T}_{\left[ P_{i}P_{i+1}%
\right] }$. The skeleton points $P_{0},P_{1}$ do not belong to the surface (%
\ref{b6.3}), which describes the interval $\mathcal{T}_{\left(
P_{i}P_{i+1}\right) }=\mathcal{T}_{\left[ P_{i}P_{i+1}\right] }\backslash
\left\{ P_{0},P_{1}\right\} $. Usually the particle mass $m$ is expressed in
the units of mass (g), and there is an universal transferring coefficient $b$%
, connecting the usual mass $m$ with the geometric mass $\mu $ 
\begin{equation}
m=b\mu =b\left\vert \mathbf{P}_{i}\mathbf{P}_{i+1}\right\vert ,\qquad \left[
b\right] =\text{g/cm}  \label{b6.11}
\end{equation}%
The same coefficient is used for connection of the geometric vector $\mathbf{%
P}_{i}\mathbf{P}_{i+1}$ with the physical momentum 4-vector $p_{k}$ 
\begin{equation}
p_{k}=bc\left( \mathbf{P}_{i}\mathbf{P}_{i+1}\right) _{k},\qquad k=0,1,2,3
\label{b6.12}
\end{equation}%
where $c$ is the speed of the light and $\left( \mathbf{P}_{i}\mathbf{P}%
_{i+1}\right) _{k}$ are coordinates of the vector $\mathbf{P}_{i}\mathbf{P}%
_{i+1}$ in some inertial coordinate system.

Thus, conceptually the dynamics of a particle is a corollary of a
geometrical description. However, the conventional dynamic formalism for
description of the free particle motion is suited only for description of
single-variant motion. In this sense the conventional single-variant
dynamics agrees with the Minkowski space-time geometry, whereas it disagrees
with the multivariant motion in the distorted space-time $V_{\mathrm{d}}$.

From the viewpoint of the distorted geometry $\mathcal{G}_{\mathrm{d}}$ the
interval $\mathcal{T}_{\left( P_{0}P_{1}\right) }$ describes the particle,
which has the shape of a hallow 3-sphere of the radius $r$ ($r_{\mathrm{\min 
}}<r<r_{\mathrm{\max }}$) . Such a spherical particle exists at rest in the
coordinate system $K_{0}$ during the proper time $t$ ($0<t<\mu /c)$. At the
proper time $t=\mu /c$ the spherical particle degenerates into the pointlike
particle, located at the point $P_{1}=\left\{ \mu ,0,0,0\right\} $. At the
next time moment $t>\mu /c$ the pointlike particle turns into spherical
particle of the radius $r$, ($r_{\mathrm{\min }}<r<r_{\mathrm{\max }}$) and
moves in the random spatial direction with the speed $\left\vert \mathbf{v}%
\right\vert =c\cdot $arth$\left( \theta _{\mathrm{M}}\right) =c\cdot $arth$%
\left( \sqrt{2}\frac{\lambda _{0}}{\mu }\right) $ in the coordinate system $%
K_{0}$. Simultaneously this spherical particle is at rest in some coordinate
system $K_{1}$, moving with respect to the coordinate system $K_{0}$ with
the velocity $\mathbf{v,}$ $\left( \left\vert \mathbf{v}\right\vert =c\cdot 
\text{arth}\left( \theta _{\mathrm{M}}\right) \right) $. The spherical
particle is at rest in the coordinate system $K_{1}$ during the proper time $%
t$, $\mu /c<t<2\mu /c$. At the proper time $t=2\mu /c$ the spherical
particle degenerates into the pointlike particle at the point $P_{2}$ and so
on. The classical mechanics cannot describe such a geometrical object as the
spherical particle without description of internal degrees of freedom. It
cannot also describe transformation of a spherical particle into a pointlike
particle and vice versa.

There is a problem of construction of the multivariant dynamics, which would
agree with the distorted geometry $\mathcal{G}_{\mathrm{d}}$, described by
the world function (\ref{a3.1}). To describe multivariant motion of a
particle, one needs to consider all variants of motion simultaneously and to
obtain some average description (some average world lines of a particle).
Such a description we shall produce on the Minkowski manifold, where results
of all mathematical operations (equality, summation, multiplication,
differentiation, etc.) are defined uniquely. Working on the Minkowski
manifold and using single-valued operation, defined on this manifold, we
shall use geometry and properties of geometrical objects, defined by the
world function (\ref{a3.1}).

In general, the classical dynamics has some experience of multivariant
motion description. In the case, when there are dynamic equations for a
single dynamic system, but there are different variants of the initial
conditions, the motion appears to be multivariant, because of different
variants of initial conditions. In this case one uses the statistical
ensemble as a dynamic system, consisting of many identical independent
dynamic systems.

For instance, a free nonrelativistic particle is described by the action%
\begin{equation}
\mathcal{A}\left[ \mathbf{x}\right] =\int m\left( \frac{d\mathbf{x}}{dt}%
\right) ^{2}dt  \label{a6.14}
\end{equation}%
where $\mathbf{x=x}\left( t\right) $ describes the world line of the
classical pointlike particle.

The statistical ensemble of free nonrelativistic particles is the dynamic
system, described by the action 
\begin{equation}
\mathcal{A}\left[ \mathbf{x}\right] =\int m\left( \frac{d\mathbf{x}}{dt}%
\right) ^{2}dtd\mathbf{\xi }  \label{a6.15}
\end{equation}%
where $\mathbf{x=x}\left( t,\mathbf{\xi }\right) $, $\mathbf{\xi }=\left\{
\xi _{1},\xi _{2},...\xi _{n}\right\} $. The independent variables $\mathbf{%
\xi }$ label particles of the statistical ensemble. The function $\mathbf{x=x%
}\left( t,\mathbf{\xi }\right) $ at fixed $\mathbf{\xi }$ describes the
world line of a particle, labelled by the label $\mathbf{\xi }$. The number $%
n$ of variables $\mathbf{\xi }$ may be arbitrary, because it is of no
importance, how the elements of the statistical ensemble are labelled.
However, usually the number of variables $\mathbf{\xi }$ is chosen to be
equal to the number of variables $\mathbf{x}=\left\{
x^{1},x^{2},x^{3}\right\} $, in order it be possible to resolve the
equations 
\begin{equation}
\mathbf{x}=\mathbf{x}\left( t,\mathbf{\xi }\right)  \label{a6.16}
\end{equation}%
in the form $\mathbf{\xi }=\mathbf{\xi }\left( t,\mathbf{x}\right) $ and to
use variables $t,\mathbf{x}$ as independent variables (Euler coordinates).

Actions (\ref{a6.14}) and (\ref{a6.16}) generate the same dynamic equations%
\begin{equation}
m\frac{d^{2}\mathbf{x}}{dt^{2}}=0  \label{a6.17}
\end{equation}%
Two dynamic systems (\ref{a6.14}) and (\ref{a6.15}) distinguish only in the
fact that the dynamic system (\ref{a6.15}) realizes a single-variant
description, whereas the statistical ensemble (\ref{a6.14}) realizes a
multivariant description, where initial conditions for different elements of
the statistical ensemble are different. Dynamic equations (\ref{a6.17}) form
a system of \textit{ordinary differential equations}.

If the motion of a single particle $\mathcal{S}_{\mathrm{st}}$ is
multivariant (stochastic), the dynamic equations for the statistical
ensemble $\mathcal{E}\left[ \mathcal{S}_{\mathrm{st}}\right] $ cease to be
ordinary differential equations. They become to be partial differential
equations, which cannot be reduced to the system of ordinary differential
equations. In this case one cannot obtain dynamic equations for a single
particle $\mathcal{S}_{\mathrm{st}}$, although there are dynamic equations
for the statistical ensemble $\mathcal{E}\left[ \mathcal{S}_{\mathrm{st}}%
\right] $.

For instance, let us consider the action of the form%
\begin{equation}
\mathcal{A}_{\mathcal{E}\left[ \mathcal{S}_{\mathrm{st}}\right] }\left[ 
\mathbf{x},\mathbf{u}\right] =\int \dint\limits_{V_{\xi }}\left\{ \frac{m}{2}%
\mathbf{\dot{x}}^{2}+\frac{m}{2}\mathbf{u}^{2}-\frac{\hbar }{2}\mathbf{%
\nabla u}\right\} dtd\mathbf{\xi },\qquad \mathbf{\dot{x}\equiv }\frac{d%
\mathbf{x}}{dt}  \label{d1.5}
\end{equation}%
The dependent variable $\mathbf{x}=\mathbf{x}\left( t,\mathbf{\xi }\right) $
describes the regular component of the particle motion. The dependent
variable $\mathbf{u}=\mathbf{u}\left( t,\mathbf{x}\right) $ describes the
mean value of the stochastic velocity component, $\hbar $ is the quantum
constant. Operator 
\begin{equation}
\mathbf{\nabla }=\left\{ \frac{\partial }{\partial x^{1}},\frac{\partial }{%
\partial x^{2}},\frac{\partial }{\partial x^{3}}\right\}  \label{d1.5a}
\end{equation}
is operator in the space of coordinates $\mathbf{x}=\left\{
x^{1},x^{2},x^{3}\right\} $.

To obtain the action functional for $\mathcal{S}_{\mathrm{st}}$ from the
action (\ref{d1.5}) for $\mathcal{E}\left[ \mathcal{S}_{\mathrm{st}}\right] $%
, we should omit integration over $\mathbf{\xi }$ in (\ref{d1.5}), as it
follows from comparison of (\ref{a6.14}) and (\ref{a6.15}). We obtain%
\begin{equation}
\mathcal{A}_{\mathcal{S}_{\mathrm{st}}}\left[ \mathbf{x},\mathbf{u}\right]
=\int \left\{ \frac{m}{2}\mathbf{\dot{x}}^{2}+\frac{m}{2}\mathbf{u}^{2}-%
\frac{\hbar }{2}\mathbf{\nabla u}\right\} dt,\qquad \mathbf{\dot{x}\equiv }%
\frac{d\mathbf{x}}{dt}  \label{d1.6}
\end{equation}%
where $\mathbf{x}=\mathbf{x}\left( t\right) $ and $\mathbf{u}=\mathbf{u}%
\left( t,\mathbf{x}\right) $ are dependent dynamic variables. The action
functional (\ref{d1.6}) is not well defined (for $\hbar \neq 0$), because
the operator $\mathbf{\nabla }$ is defined in some 3-dimensional vicinity of
point $\mathbf{x}$, but not at the point $\mathbf{x}$ itself. As far as the
action functional (\ref{d1.6}) is not well defined, one cannot obtain
dynamic equations for $\mathcal{S}_{\mathrm{st}}$. By definition it means
that the motion of the particle $\mathcal{S}_{\mathrm{st}}$ is multivariant
(stochastic). Setting $\hbar =0$ in (\ref{d1.6}), we transform the action (%
\ref{d1.6}) into the action (\ref{a6.14}), because in this case $\mathbf{u}%
=0 $ in virtue of dynamic equations.

Let us return to the action (\ref{d1.5}) and obtain dynamic equations for
the statistical ensemble $\mathcal{E}\left[ \mathcal{S}_{\mathrm{st}}\right] 
$ of physical systems $\mathcal{S}_{\mathrm{st}}$. Variation of (\ref{d1.5})
with respect to $\mathbf{u}$ gives%
\begin{eqnarray*}
\delta \mathcal{A}_{\mathcal{E}\left[ \mathcal{S}_{\mathrm{st}}\right] }%
\left[ \mathbf{x},\mathbf{u}\right] &=&\int \dint\limits_{V_{\xi }}\left\{ m%
\mathbf{u}\delta \mathbf{u}-\frac{\hbar }{2}\mathbf{\nabla }\delta \mathbf{u}%
\right\} dtd\mathbf{\xi } \\
&=&\int \dint\limits_{V_{\mathbf{x}}}\left\{ m\mathbf{u}\delta \mathbf{u}-%
\frac{\hbar }{2}\mathbf{\nabla }\delta \mathbf{u}\right\} \frac{\partial
\left( \xi _{1},\xi _{2},\xi _{3}\right) }{\partial \left(
x^{1},x^{2},x^{3}\right) }dtd\mathbf{x} \\
&=&\int \dint\limits_{V_{\mathbf{x}}}\delta \mathbf{u}\left\{ m\mathbf{u}%
\rho +\frac{\hbar }{2}\mathbf{\nabla }\rho \right\} dtd\mathbf{x-}\int
\doint \frac{\hbar }{2}\rho \delta \mathbf{u}dtd\mathbf{S}
\end{eqnarray*}%
where%
\begin{equation}
\rho =\frac{\partial \left( \xi _{1},\xi _{2},\xi _{3}\right) }{\partial
\left( x^{1},x^{2},x^{3}\right) }=\left( \frac{\partial \left(
x^{1},x^{2},x^{3}\right) }{\partial \left( \xi _{1},\xi _{2},\xi _{3}\right) 
}\right) ^{-1}  \label{d1.7a}
\end{equation}%
We obtain the following dynamic equation%
\begin{equation}
m\rho \mathbf{u}+\frac{\hbar }{2}\mathbf{\nabla }\rho =0,  \label{d1.7}
\end{equation}%
Variation of (\ref{d1.5}) with respect to $\mathbf{x}$ gives%
\begin{equation}
m\frac{d^{2}\mathbf{x}}{dt^{2}}=\mathbf{\nabla }\left( \frac{m}{2}\mathbf{u}%
^{2}-\frac{\hbar }{2}\mathbf{\nabla u}\right)  \label{d1.8}
\end{equation}%
Here $d/dt$ means the substantial derivative with respect to time $t$%
\begin{equation*}
\frac{dF}{dt}\equiv \frac{\partial \left( F,\xi _{1},\xi _{2},\xi
_{3}\right) }{\partial \left( t,\xi _{1},\xi _{2},\xi _{3}\right) }
\end{equation*}

Resolving (\ref{d1.7}) with respect to $\mathbf{u}$, we obtain the equation%
\begin{equation}
\mathbf{u}=-\frac{\hbar }{2m}\mathbf{\nabla }\ln \rho ,  \label{d1.9}
\end{equation}%
which reminds the expression for the mean velocity of the Brownian particle
with the diffusion coefficient $D=\hbar /2m$.

Eliminating the velocity $\mathbf{u}$ from dynamic equations (\ref{d1.8})
and (\ref{d1.9}) and going to independent Eulerian variables $t,\mathbf{x}$,
we obtain the dynamic equations of the hydrodynamic type for the mean motion
of the stochastic particle $\mathcal{S}_{\mathrm{st}}$%
\begin{equation}
m\frac{d^{2}\mathbf{x}}{dt^{2}}=m\left( \frac{\partial \mathbf{v}}{\partial t%
}+\left( \mathbf{v\nabla }\right) \mathbf{v}\right) =-\mathbf{\nabla }U_{%
\mathrm{B}},\qquad \frac{\partial \rho }{\partial t}+\mathbf{\nabla }\left(
\rho \mathbf{v}\right) =0  \label{d1.10}
\end{equation}%
where $\mathbf{v=}\frac{d\mathbf{x}}{dt}$ and $U_{\mathrm{B}}$ is the Bohm
potential \cite{B52}%
\begin{equation}
U_{\mathrm{B}}=U\left( \rho ,\mathbf{\nabla }\rho ,\mathbf{\nabla }^{2}\rho
\right) =\frac{\hbar ^{2}}{8m}\frac{\left( \mathbf{\nabla }\rho \right) ^{2}%
}{\rho ^{2}}-\frac{\hbar ^{2}}{4m}\frac{\mathbf{\nabla }^{2}\rho }{\rho }
\label{d1.11}
\end{equation}

In the case of the irrotational flow the hydrodynamic equations (\ref{d1.10}%
) are equivalent to the Schr\"{o}dinger equations \cite{B52}. Thus, the
action (\ref{d1.5}) and dynamic equations (\ref{d1.10}) describe correctly
multivariant motion of particles of statistical ensemble $\mathcal{E}\left[ 
\mathcal{S}_{\mathrm{st}}\right] $. But it remains unclear, how the form of
the action (\ref{d1.5}) is connected with the multivariant space-time
geometry. This question was investigated in \cite{R91}. It was shown that
for agreement of the space-time geometry with the action (\ref{d1.5}) the
distortion $d=\lambda _{0}^{2}$ in the world function (\ref{a1.3}) should be
chosen in the form 
\begin{equation}
d=\lambda _{0}^{2}=\frac{\hslash }{2bc}  \label{d1.12}
\end{equation}%
where $\hbar $ is the quantum constant, $c$ is the speed of the light and $b$
is the constant from the relation (\ref{b6.11}), connecting the geometric
mass $\mu $ with the usual mass $m$ of the particle.

Here we suggest only some simple arguments for explanation of the connection
of the world function of the space-time and the action (\ref{d1.5}),
describing motion the ensemble of free particles.

As it follows from the relation (\ref{b6.10}), the particle velocity,
described by the link $\mathcal{T}_{\left[ P_{1}P_{2}\right] }$, has two
components. One component $\mathbf{v}_{\mathrm{reg}}=\mathbf{p}/m$ is
regular. It is determined by the particle momentum $\mathbf{p}$. Other
component of velocity $\mathbf{v}_{\mathrm{st}}$ is conditioned by the
random walk of the particle. Its average value $\left\langle \mathbf{v}_{%
\mathrm{st}}\right\rangle $ depends on the state of the whole ensemble. It
has the form 
\begin{equation}
\left\langle \mathbf{v}_{\mathrm{st}}\right\rangle =-\alpha r_{\min }c\theta
_{\mathrm{M}}\mathbf{\nabla }\log \rho ,\qquad r_{\mathrm{\min }}=\sqrt{%
\frac{3}{2}}\lambda _{0},\qquad \theta _{\mathrm{M}}=\sqrt{2}\frac{\lambda
_{0}}{\mu }  \label{a6.18}
\end{equation}%
where $\alpha $ is some real number of the order of $1$ and $\rho $ is the
density of particles in the statistical ensemble.

On one hand, it follows from the action (\ref{d1.5}), which gives the true
description of the mean particle motion, that the mean stochastic velocity
can be presented in the form (\ref{d1.9}). On the other hand, comparison of
relations (\ref{a6.18}) and (\ref{d1.9}) leads to result (\ref{d1.12}). It
is important, that constant $b$ does not appear in the action (\ref{d1.5}),
and one cannot determine the value of $b$ as well the value of $q=\lambda
_{0}^{2}$ experimentally. Thus, for explanation of quantum effects it is
important the fact of the multivariance existence, but not its numerical
value.

\section{Concluding remarks}

The non-Euclidean method of the generalized geometry construction admits one
to construct all possible generalized geometries. These geometries may be
continuous or discrete, They may have alternating dimension, or have no
dimension at all. The non-Euclidean method deals only with the world
function, which is the only essential characteristic of geometry. The
non-Euclidean method does not use coordinate description, and there is no
necessity to take into account and remove arbitrariness, connected with a
usage of coordinates.

The non-Euclidean method does not need a separation of geometrical
propositions into axioms and theorems. It does not need a test of the axioms
consistency, connected with this separation. The non-Euclidean method admits
one to discover the property of multivariance, which is very general
property of generalized geometries. The multivariance appears to be a very
important property of the space-time geometry, responsible for quantum
effects. Existence of multivariance dictates a new revision of the
space-time geometry. The multivariance of the space-time geometry admits one
to move along the way of the further physics geometrization. The particle
mass appears to be a geometrical characteristic of a particle. It becomes to
put the question on existence of geometrical objects. In particular, it
becomes to be possible to consider the confinement problem as the
geometrical problem of the complicated object existence, but not as a
dynamical problem of the pointlike particles confinement inside a restricted
volume.

\end{document}